# TESTING CONDITIONAL INDEPENDENCE VIA ROSENBLATT TRANSFORMS

BY KYUNGCHUL SONG

*University of Pennsylvania*

This paper proposes new tests of conditional independence of two random variables given a single-index involving an unknown finite-dimensional parameter. The tests employ Rosenblatt transforms and are shown to be distribution-free while retaining computational convenience. Some results from Monte Carlo simulations are presented and discussed.

**1. Introduction.** Suppose that $Y$ and $Z$ are random variables, and let $\lambda_\theta(X)$ be a real function of a random vector $X$ indexed by a parameter $\theta \in \Theta \subset \mathbf{R}^d$. The function $\lambda_\theta(\cdot)$ is known up to $\theta \in \Theta$. For example, we may consider $\lambda_\theta(X) = h(X^\top \theta)$ for some known function $h$. Suppose that an estimable parameter $\theta_0 \in \Theta$ is given. This paper proposes a distribution-free method of testing conditional independence of $Y$ and $Z$ given $\lambda_{\theta_0}(X)$,

$$(1) \qquad Y \perp\!\!\!\perp Z | \lambda_{\theta_0}(X).$$

When $Y$ and $Z$ are conditionally independent given $\lambda_{\theta_0}(X)$, it means that "learning the value of $Z$ does not provide additional information about $Y$, once we know $\lambda_{\theta_0}(X)$" [Pearl (2000), page 11]. Hence conditional independence is a central notion in modeling causal relations, and its importance in graphical modeling is widely known [e.g., Lauritzen (1996), Pearl (2000)]. In the literature of program evaluations, testing conditional independence of the observed outcome and the treatment decision given observed covariates can serve as testing lack of treatment effects under the assumption of strong ignorability [Heckman, Ichimura and Todd (1997)]. Conditional independence is, sometimes, a direct implication of economic theory. For example, in the literature of insurance, the presence of positive conditional dependence between coverage and risk is known to be a direct consequence of

---

adverse selection under information asymmetry [e.g., Chiappori and Salanié (2000)].

Testing independence for continuous variables has drawn the attention of many researchers. To name but a few, see Hoeffding (1948), Blum, Kiefer and Rosenblatt (1961), Skaug and Tjøstheim (1993), Robinson (1991), Delgado and Mora (2000) and Hong and White (2005). There has also been a growth of interest in testing extremal dependence. See a recent paper by Zhang (2008) and references therein. In contrast, the literature of testing conditional independence for continuous variables appears rather recent and includes relatively few researches. See Linton and Gozalo (1997), Delgado and González Manteiga (2001), Angrist and Kuersteiner (2004) and Su and White (2008), among others. None of these tests focuses on conditional independence between continuous variables with unknown $\theta_0$ and is distribution-free at the same time.

Distribution-free tests have asymptotic critical values that do not change as we move from one probability to another within the null hypothesis. Many goodness-of-fit tests that have nontrivial asymptotic power against $\sqrt{n}$-converging Pitman local alternatives are not distribution-free. To deal with this problem, the literature either suggests the use of approximate critical values through bootstrap or the transformation of the test using the innovation martingale approach pioneered by Khmaladze (1993).

This paper shows that for testing conditional independence, we can generate distribution-free tests by appropriately using Rosenblatt transforms—a multivariate version of a probability integral transform studied by Rosenblatt (1952). Based on the result, this paper proposes a bootstrap method that is computationally attractive. This bootstrap procedure does not require the re-estimation of $\theta_0$ for each bootstrap sample. This is convenient when the dimension of $\theta_0$ is large and its estimation involves numerical optimization.

The Rosenblatt transform is closely related to the probability integral transform of a single-index suggested by Stute and Zhu (2005). However, the nature of the problem is distinguished from that of Stute and Zhu (2005). First, our test of conditional independence is both omnibus (when the test is two-sided) and distribution-free, while the single-index restriction test of Stute and Zhu (2005) fails to be distribution-free when it is designed to be omnibus. This is purely due to the nature of conditional independence as distinct from a single-index restriction. Second, our test contains the probability integral transform inside functions that are potentially discontinuous, making it cumbersome to rely on the $U$-process theory [e.g., de la Peña and Giné (1999)]. This paper deals with this difficulty by directly establishing the bracketing entropy bounds for functions involving the probability integral transforms. These entropy bounds can be used for other purposes.

This paper is organized as follows. In the next section, we introduce the basic assumptions and test statistics, and develop asymptotic theory. In Section 3, we consider a bootstrap method. Section 4 deals with the case where $Z$ is discrete. In Section 5, we present and discuss the results from the Monte Carlo simulation study. In the Appendix, we offer the mathematical proofs.

**2. Main results.** Suppose that we are given a random vector $(Y, Z, X)$ distributed by $P$ and a real valued function $\lambda_{\theta_0}(\cdot)$ on $\mathbf{R}^{d_X}$ which is known up to a parameter $\theta_0 \in \Theta \subset \mathbf{R}^d$. For brevity, let $\lambda_0(\cdot) \equiv \lambda_{\theta_0}(\cdot)$. Throughout this paper, we assume that $\lambda_0(X)$ is continuous. Define $F_0(\cdot)$ to be the distribution function of $\lambda_0(X)$ and let $U \equiv F_0(\lambda_0(X))$, $F_{Y|U}(\cdot|U) \equiv P\{Y \leq \cdot|U\}$, and $F_{Z|U}(\cdot|U) \equiv P\{Z \leq \cdot|U\}$. Then, the main focus of this paper is on testing the following hypothesis:

$$(2) \qquad H_0 : P\{Y \leq y, Z \leq z|U\} = F_{Y|U}(y|U)F_{Z|U}(z|U) \qquad \text{wp } 1,$$

for all $(y, z)$ in the support of $(Y, Z)$. The notation "wp 1" means that the statement holds with probability one with respect to the distribution of $U$. Certainly, this hypothesis is equivalent to (1) with probability one because $F_0(\cdot)$ is continuous.

Throughout the paper, the norm $\|\cdot\|_\infty$ represents the sup norm, $\|\cdot\|$, the Euclidean norm and $\|\cdot\|_{P,p}$, the $L_p(P)$-norm. Let $B(\theta_0, \delta) \equiv \{\theta \in \Theta : \|\theta - \theta_0\| \leq \delta\}$. Define

$$\tilde{Y} \equiv F_{Y|U}(Y|U) \quad \text{and} \quad \tilde{Z} \equiv F_{Z|U}(Z|U),$$

then $(\tilde{Z}, U)$ is distributed as the joint distribution of two independent uniform $[0, 1]$ random variables, and so is $(\tilde{Y}, U)$, if $(Z, \lambda_0(X))$ and $(Y, \lambda_0(X))$ are continuous [Rosenblatt (1952)]. The transform of $(Z, \lambda_0(X))$ into $(\tilde{Z}, U)$ is called *the Rosenblatt transform*, due to Rosenblatt (1952). Let $f_{\tilde{Y}|Z,\theta}(y|z, \bar{\lambda}_1, \bar{\lambda}_0)$ be the conditional density of $\tilde{Y}$ given $(\tilde{Z}, \lambda_\theta(X), \lambda_0(X)) = (z, \bar{\lambda}_1, \bar{\lambda}_0)$ with respect to a $\sigma$-finite conditional measure. We also define $f_{\tilde{Z}|Y,\theta}(z|y, \bar{\lambda}_1, \bar{\lambda}_0)$ similarly by interchanging the roles of $\tilde{Y}$ and $\tilde{Z}$.

ASSUMPTION 1. (i) $(Y, Z, \lambda_0(X))$ is continuous.

(ii) For some $\delta > 0$,

(a) $\lambda_\theta(\cdot)$, $\theta \in B(\theta_0, \delta)$, is uniformly bounded and Lipschitz in $\theta$, that is, for any $\theta_1, \theta_2 \in B(\theta_0, \delta)$,

$$\|\lambda_{\theta_1} - \lambda_{\theta_2}\|_\infty \leq C\|\theta_1 - \theta_2\| \qquad \text{for some } C > 0 \quad \text{and}$$

(b) $\lambda_\theta(X)$ is continuous, having a density function bounded uniformly over $\theta \in B(\theta_0, \delta)$.

(iii) For some $\delta > 0$, $f_{Z|Y,\theta}(z|y, \cdot, \bar{\lambda}_0)$ and $f_{Y|Z,\theta}(y|z, \cdot, \bar{\lambda}_0)$ are continuously differentiable with derivatives bounded uniformly over $(y, z, \bar{\lambda}_0, \theta) \in [0,1]^2 \times \mathbf{R} \times B(\theta_0, \delta)$.

Later in the paper (Section 4), we deal with the case where either $Y$ or $Z$ is discrete. The uniform boundedness condition in (ii) is innocuous, because by choosing a strictly increasing function $\Phi$ on $[0,1]$, we can redefine $\lambda'_\theta = \Phi \circ \lambda_\theta$. The Lipschitz continuity in $\theta$ can be made to hold by choosing this $\Phi$ appropriately. The absolute continuity condition in (ii)(b) is satisfied in particular when $\lambda_\theta(X) = h(X^\top \theta)$ with a continuous, strictly increasing function $h$ and $X^\top \theta$ is continuous.

Define
$$\gamma_z(\cdot) \equiv z \exp(\cdot \times z) \quad \text{and} \quad \gamma_z^\perp(\cdot) \equiv \gamma_z(\cdot) - \{\exp(z) - 1\}.$$

For a class of functions $\beta_u(\cdot)$, $u \in [0,1]$, consider the following null hypothesis:

(3) $$H_0 : \mathbf{E}[\beta_u(U)\gamma_z^\perp(\tilde{Z})\gamma_y^\perp(\tilde{Y})] = 0 \qquad \forall (u, y, z) \in [0,1]^3.$$

The lemma below establishes that under Assumption 1(i), and an appropriate condition for $\beta_u$, the null hypothesis in (2) and the null hypothesis in (3) are equivalent. The result relies on Lemma 1 of Bierens (1990).

LEMMA 1. *Suppose that Assumption 1(i) is satisfied. Furthermore, assume that the class $\{\beta_u, u \in [0,1]\}$ is such that (3) implies the following:*

(4) $$\mathbf{E}[\gamma_z^\perp(\tilde{Z})\gamma_y^\perp(\tilde{Y})|U] = 0 \qquad \text{wp } 1 \ \forall (y, z) \in [0,1]^2.$$

*Then, the hypothesis in (2) and the hypothesis in (3) are equivalent.*

PROOF. It is easy to see that the conditional independence (2) implies (3). We prove the converse. First, we show that the conditional independence of $\tilde{Y}$ and $\tilde{Z}$ given $U$ implies (2). Suppose that this conditional independence holds. Let $\tilde{y} = F_{Y|U}(y|U)$ and $\tilde{z} = F_{Z|U}(z|U)$ for brevity. Write

(5)
$$\begin{aligned}
P\{&\tilde{Y} \leq \tilde{y}, \tilde{Z} \leq \tilde{z}|U\} \\
&= P\{\tilde{Y} \leq \tilde{y}, \tilde{Z} \leq \tilde{z}, Y \leq y, Z \leq z|U\} \\
&\quad + P\{\tilde{Y} \leq \tilde{y}, \tilde{Z} \leq \tilde{z}, Y > y, Z \leq z|U\} \\
&\quad + P\{\tilde{Y} \leq \tilde{y}, \tilde{Z} \leq \tilde{z}, Y \leq y, Z > z|U\} \\
&\quad + P\{\tilde{Y} \leq \tilde{y}, \tilde{Z} \leq \tilde{z}, Y > y, Z > z|U\}.
\end{aligned}$$

Following Angus (1994), the second probability on the right-hand side is bounded by

$$\begin{aligned}
P\{\tilde{Y} \leq \tilde{y}, Y > y|U\} &= P\{\tilde{Y} = \tilde{y}, Y > y|U\} \\
&= P\{F_{Y|U}(Y|U) = F_{Y|U}(y|U), Y > y|U\} = 0,
\end{aligned}$$

because conditional on $U = u$, the event in the last probability is contained in the event of $Y$ lying in the interior of an interval of constancy of $F_{Y|U}(\cdot|u)$. The conditional probability measure of this event is certainly zero. Similarly, the last two probabilities in (5) can also be shown to be zero. If $Y \leq y$ and $Z \leq z$, then $\tilde{Y} \leq \tilde{y}$ and $\tilde{Z} \leq \tilde{z}$. Therefore, we obtain from (5) that

$$(6) \qquad P\{Y \leq y, Z \leq z|U\} = P\{\tilde{Y} \leq \tilde{y}, \tilde{Z} \leq \tilde{z}|U\}.$$

Using a similar argument, we can also obtain that

$$(7) \qquad \begin{aligned} P\{Y \leq y|U\} &= P\{\tilde{Y} \leq \tilde{y}|U\} = \tilde{y} \quad \text{and} \\ P\{Z \leq z|U\} &= P\{\tilde{Z} \leq \tilde{z}|U\} = \tilde{z}, \end{aligned}$$

because $\tilde{Y}$ is uniformly distributed on $[0,1]$ and independent of $U$, and so is $\tilde{Z}$. Conditional independence of $\tilde{Y}$ and $\tilde{Z}$ given $U$ implies (2) through (6) and (7).

Now, we show that (3) implies conditional independence of $\tilde{Y}$ and $\tilde{Z}$ given $U$. Let $f(t_1, t_2|u) : [0,1]^2 \to [0, \infty)$ be the conditional density of $(\tilde{Z}, \tilde{Y})$ given $U = u$. Through (4), (3) implies that for all $(z, y) \in [0,1]^2$,

$$zy \int_0^1 \int_0^1 e^{zt_1 + yt_2} \{f(t_1, t_2|U) - 1\} \, dt_1 \, dt_2 = 0 \qquad \text{wp } 1.$$

By Lemma 1 of Bierens (1990) [see also Stinchcombe and White (1998), page 4], $f(t_1, t_2|U) = 1\{(t_1, t_2) \in [0,1]^2\}$, a.e., for almost every $(t_1, t_2) \in [0,1]^2$, yielding conditional independence of $\tilde{Y}$ and $\tilde{Z}$ given $U$. □

The condition for $\beta_u(\cdot)$ in Lemma 1 is explained in Stinchcombe and White (1998). For example, the choice of $\beta_u(U) = 1\{U \leq u\}$ or $\beta_u(U) = \exp(Uu)$ satisfies this condition [see Bierens (1990), Lemma 1, for the latter choice]. From now on, we assume that $\beta_u(\cdot)$ satisfies the condition in Lemma 1 and focus on the null hypothesis in (3). This condition for $\beta_u(\cdot)$ is not used for the weak convergence theory in Theorem 1 below.

ASSUMPTION 2. (i) $\beta_u(\cdot), u \in [0,1]$, is uniformly bounded in $[0,1]$, and for each $u \in [0,1], \beta_u(\cdot)$ is of bounded variation.

(ii) $F_{Y|U}(y|\cdot)$ and $F_{Z|U}(z|\cdot)$ are twice continuously differentiable with derivatives bounded uniformly over $(z, y) \in \mathbf{R}^2$.

Assumption 2(i) is very weak and satisfied by most functions used in the literature. This flexibility in choosing the class $\beta_u$ is important because the choice of $\beta_u$ plays a significant role in determining the asymptotic power properties of the test in general. Assumption 2(ii) is analogous to Condition A(i) in Theorem 2.1 of Stute and Zhu (2005) or A2 of Delgado and González Manteiga (2001) on page 1475.

Throughout this paper, we assume that the observations $(Y_i, X_i, Z_i)_{i=1}^n$ are i.i.d. from $P$. We also assume that the parameter $\theta_0$ is identified from data and estimable. For example, in the literature of program evaluations, this assumption is satisfied because the parameter $\theta_0$ constitutes the single-index in the propensity score. Let $\hat{\theta}$ be a consistent estimator of $\theta_0$, and define

$$\hat{U}_i \equiv F_{n,\hat{\theta},i}(\lambda_{\hat{\theta}}(X_i)), \qquad \hat{Z}_i \equiv \hat{F}_{Z|U,i}(Z_i|\hat{U}_i) \quad \text{and} \quad \hat{Y}_i \equiv \hat{F}_{Y|U,i}(Y_i|\hat{U}_i),$$

where $F_{n,\hat{\theta},i}(\bar{\lambda}) \equiv \frac{1}{n-1}\sum_{j=1,j\neq i}^n 1\{\lambda_{\hat{\theta}}(X_j) \leq \bar{\lambda}\}$, and

$$(8) \qquad \hat{F}_{Y|U,i}(y|u) \equiv \frac{\sum_{j=1,j\neq i}^n 1\{Y_j \leq y\}K_h(\hat{U}_j - u)}{\sum_{j=1,j\neq i}^n K_h(\hat{U}_j - u)},$$

where $K_h(x) = K(x/h)/h$, $K(\cdot)$ is a kernel function and $h$ is the bandwidth parameter. We similarly define $\hat{F}_{Z|U,i}(z|u)$. As for the estimator $\hat{\theta}$, the kernel and the bandwidth, we assume the following.

ASSUMPTION 3. (i) $\|\hat{\theta} - \theta_0\| = O_P(n^{-1/2})$.

(ii) (a) $K$ is symmetric, nonnegative, twice continuously differentiable, has a compact support, and $\int_{-\infty}^{\infty} K(s) = 1$.

(b) $h = Cn^{-s}$ with $1/6 < s < 1/4$ for some $C > 0$.

When $\hat{\theta}$ is an $M$-estimator, the rate of convergence in (i) can be obtained following the procedure of Theorem 3.2.5 of van der Vaart and Wellner (1996). The estimation method of $\theta_0$ depends on a further specification of the testing environment. For example, the conditioning variable $\lambda_\theta(X_i)$ may originate from the nonlinear regression model,

$$W_i = \lambda_{\theta_0}(X_i) + \varepsilon_i,$$

where $\varepsilon_i$ satisfies $\mathbf{E}[\varepsilon_i|X_i] = 0$. The $\sqrt{n}$-consistent estimation of $\theta_0$ in this case is well known in the literature [see, e.g., van de Geer (2000)]. Assumption 3(ii)(a) is used by Stute and Zhu (2005). Unlike their procedure, the bandwidth condition in (b) does not require undersmoothing.

Define the infeasible and feasible processes

$$\nu_n(r) \equiv \frac{1}{\sqrt{n}}\sum_{i=1}^n \beta_u(U_i)\gamma_z^\perp(\tilde{Z}_i)\gamma_y^\perp(\tilde{Y}_i)$$

and

$$\hat{\nu}_n(r) \equiv \frac{1}{\sqrt{n}}\sum_{i=1}^n \beta_u(\hat{U}_i)\gamma_z^\perp(\hat{Z}_i)\gamma_y^\perp(\hat{Y}_i).$$

In the following, we establish weak convergence of both processes. The main complication is that the condition for $\beta_u$ (Assumption 2) is too weak to resort to linearization in handling the estimation error of $\hat{U}_i$. We deal with this difficulty by establishing a bracketing entropy bound for functions composite with bounded variation functions (see Lemma A1 in the Appendix). Let $l_\infty([0,1]^3)$ denote the space of real functions on $[0,1]^3$ that are bounded, and endowed with a sup norm $\|\cdot\|_\infty$ defined by $\|f\|_\infty = \sup_{u\in[0,1]^3}|f(u)|$. The notation $\rightsquigarrow$ denotes weak convergence in $l_\infty([0,1]^3)$ in the sense of Hoffman–Jorgensen [e.g., van der Vaart and Wellner (1996)]. Let $\langle\cdot,\cdot\rangle$ be the inner product on $L_2(du) \times L_2(du)$ defined by $\langle f,g\rangle = \int_0^1 f(u)g(u)\,du$.

THEOREM 1. *Suppose that Assumptions 1–3 hold. Then the following holds:*

(i) $\sup_{r\in[0,1]^3}|\hat{\nu}_n(r) - \nu_n(r)| = o_P(1)$, *both under $H_0$ in (3) and under Pitman local alternatives $\tilde{P}_n$ such that for some functions $a_j : [0,1]^3 \rightarrow [-1,1]$, $j = 1, 2$,*

$$\mathbf{E}_{\tilde{P}_n}[\gamma_z^\perp(\tilde{Z})|\tilde{Y}=y, U=u] = n^{-1/2}a_1(z,y,u)$$

*and*

$$\mathbf{E}_{\tilde{P}_n}[\gamma_y^\perp(\tilde{Y})|\tilde{Z}=z, U=u] = n^{-1/2}a_2(z,y,u),$$

*where $\mathbf{E}_{\tilde{P}_n}[\cdot|\tilde{Y}=y, U=u]$ and $\mathbf{E}_{\tilde{P}_n}[\cdot|\tilde{Z}=z, U=u]$ denote conditional expectations under $\tilde{P}_n$.*

(ii) $\hat{\nu}_n \rightsquigarrow \nu$ *in $l_\infty([0,1]^3)$, under $H_0$ in (3), where $\nu$ is a centered Gaussian process whose covariance kernel is given by*

$$c(r_1; r_2) = \langle\beta_{u_1}, \beta_{u_2}\rangle\langle\gamma_{z_1}^\perp, \gamma_{z_2}^\perp\rangle\langle\gamma_{y_1}^\perp, \gamma_{y_2}^\perp\rangle.$$

The asymptotic representation in (i) shows an interesting fact that $\hat{\nu}_n(r)$ is asymptotically equivalent to $\nu_n(r)$. Remarkably, the estimation error in $\hat{\theta}$ does not play a role in shaping the asymptotic distribution of the process $\hat{\nu}_n(r)$. This finding is analogous to what Stute and Zhu (2005) found in the context of testing a single-index restriction.

Based on the result in Theorem 1, we can construct a test statistic

$$(9) \qquad T_n = \Gamma\hat{\nu}_n$$

by taking a continuous functional $\Gamma$. For example, in the case of two sided tests, we may take

$$(10) \qquad \Gamma_{\mathrm{KS}}\hat{\nu}_n = \sup_{r\in[0,1]^3}|\hat{\nu}_n(r)| \quad \text{or} \quad \Gamma_{\mathrm{CM}}\hat{\nu}_n = \left(\int_{[0,1]^3}\hat{\nu}_n(r)^2\,dr\right)^{1/2}.$$

The first example is of Kolmogorov–Smirnov-type and the second one is of Cramér–von Mises-type. Asymptotic unbiasedness for these tests against $\sqrt{n}$-converging local alternatives can be established using Anderson's lemma. In the case of one-sided tests, we may take

$$\Gamma_{\mathrm{KS}}^{+}\hat{\nu}_n = \sup_{r\in[0,1]^3} \hat{\nu}_n(r) \quad \text{or} \quad \Gamma_{\mathrm{CM}}^{+}\hat{\nu}_n = \left(\int_{[0,1]^3} \max\{\hat{\nu}_n(r),0\}^2\, dr\right)^{1/2}.$$

The asymptotic properties of the tests follow from Theorem 1. Indeed, under $H_0$,

$$(11) \qquad\qquad T_n = \Gamma\hat{\nu}_n \to_d \Gamma\nu.$$

This test is distribution-free, as the limiting distribution of $\Gamma\nu$ does not depend on the data generating process under $H_0$.

**3. Bootstrap tests.** The tests introduced so far are distribution-free, but, in many cases, it is not known how to simulate the Gaussian process $\nu$. In this section, we suggest a wild bootstrap method in a spirit similar to Delgado and González Manteiga (2001) [see also, among others, Härdle and Mammen (1993), Stute, González Manteiga and Quindimil (1998)].

Let $(\{\omega_{i,b}\}_{i=1}^{n})_{b=1}^{B}$ be an i.i.d. sequence of random variables that are bounded, independent of $\{Y_i, Z_i, X_i\}$, $\mathbf{E}(\omega_{i,b}) = 0$ and $\mathbf{E}(\omega_{i,b}^2) = 1$. For example, one can take $\omega_{i,b}$ with a two-point distribution assigning masses $(\sqrt{5}+1)/(2\sqrt{5})$ and $(\sqrt{5}-1)/(2\sqrt{5})$ to the points $-(\sqrt{5}-1)/2$ and $(\sqrt{5}+1)/2$. Let

$$\nu_{n,b}^{*}(r) = \frac{1}{\sqrt{n}}\sum_{i=1}^{n} \omega_{i,b}\beta_u(\hat{U}_i)\gamma_z^{\perp}(\hat{Z}_i)\gamma_y^{\perp}(\hat{Y}_i), \qquad b=1,\dots,B.$$

The bootstrap empirical process $\nu_{n,b}^{*}(r)$ is similar to those proposed by Delgado and González Manteiga (2001). Given a functional $\Gamma$, we can define bootstrap test statistics $T_{n,b}^{*} = \Gamma\nu_{n,b}^{*}$, $b=1,\dots,B$. An $\alpha$-level critical value is approximated by $c_{\alpha,n,B} = \inf\{t : B^{-1}\Sigma_{b=1}^{B}\mathbf{1}\{T_{n,b}^{*} \le t\} \ge 1-\alpha\}$, yielding bootstrap test $\mathbf{1}\{T_n > c_{\alpha,n,B}\}$, where $T_n$ is as defined in (9).

Let $F_{\Gamma\nu}$ be the distribution of $\Gamma\nu$ and let $F_{T_n^{*}}^{*}$ denote the conditional distributions of bootstrap test statistics $T_n^{*}$. Define $d(\cdot,\cdot)$ to be a distance metrizing weak convergence on the real line. [For an introductory exposition about the weak convergence of bootstrap empirical processes, see Giné (1997)]. The weak convergence follows eventually as a consequence of the almost sure multiplier CLT of Ledoux and Talagrand (1988).

THEOREM 2. *Suppose that the conditions of Theorem 1 hold under $H_0$ in (3). Then under $H_0$,*

$$d(F_{T_n^{*}}^{*}, F_{\Gamma\nu}) \to 0 \qquad in\ P.$$

The wild bootstrap procedure is easy to implement. In particular, one does not need to re-estimate $\theta_0$ or $(\tilde{Z}_i, \tilde{Y}_i)$ using each bootstrap sample. It is worth noting that this desirable property is made possible by our transforming the test into a distribution-free one.

**4. Discrete random variables.** The development so far has assumed that $Y$ and $Z$ are continuous random variables. In many important applications of conditional independence, either $Y$ or $Z$ is discrete, or more often, binary. For example, in the literature of program evaluations, the conditional independence restriction involves a binary variable representing the incidence of treatment.

From now on, we assume that $Y$ is continuous and $Z$ is discrete, taking values from a known, finite set $\mathcal{Z}$. We introduce $(\tilde{Y}, U)$ as before. Define $p_z(U) = P\{Z = z | U\}$ for $z \in \mathcal{Z}$. Similarly as in Lemma 1, we can show that the null hypothesis in (2) is equivalent to

$$H_0 : \mathbf{E}[\beta_u(U)\{1\{Z = z\} - p_z(U)\}\gamma_y^\perp(\tilde{Y})] = 0 \qquad \forall (y, u, z) \in [0,1]^2 \times \mathcal{Z},$$

if $(Y, \lambda_0(X))$ is continuous and $\beta_u$ satisfies approximate conditions similarly as in Lemma 1. We substitute the following for Assumption 1.

ASSUMPTION 1D. (i) $(Y, \lambda_0(X))$ is continuous.

(ii) Assumption 1(ii) holds for $\lambda_\theta, \theta \in \Theta$.

(iii) For some $\delta > 0$, $f_{Y|Z,\theta}(y|z, \cdot, \bar{\lambda}_0)$ is continuously differentiable with a derivative bounded uniformly over $(y, z, \bar{\lambda}_0, \theta) \in [0,1] \times \mathcal{Z} \times \mathbf{R} \times B(\theta_0, \delta)$.

(iv) For some $\varepsilon > 0$, $p_z(u) \in (\varepsilon, 1 - \varepsilon)$ for all $(u, z) \in [0,1] \times \mathcal{Z}$.

Let $\hat{p}_{z,i}(u)$ be a kernel estimator of $p_z(u)$,

$$(12) \qquad \hat{p}_{z,i}(u) = \frac{\sum_{j=1, j \neq i}^n 1\{Z_j = z\} K_h(\hat{U}_j - u)}{\sum_{j=1}^n K_h(\hat{U}_j - u)}$$

and consider the following process: for $(u, y, z) \in [0,1]^2 \times \mathcal{Z}$,

$$\bar{\nu}_n(u, y, z) \equiv \frac{1}{\sqrt{n}} \sum_{i=1}^n \frac{\beta_u(\hat{U}_i)\{1\{Z_i = z\} - \hat{p}_{z,i}(\hat{U}_i)\}\gamma_y^\perp(\hat{Y}_i)}{\sqrt{\hat{p}_{z,i}(\hat{U}_i) - \hat{p}_{z,i}(\hat{U}_i)^2}},$$

where $\hat{Y}_i$ and $\hat{U}_i$ are as defined before.

THEOREM 3. *Suppose that Assumptions 1D, 2 and 3 hold. Furthermore, the conditions for $K$ and $h$ used for $\hat{p}_{z,i}(\cdot)$ are the same as Assumption 3(ii). Then under $H_0$,*

$$\bar{\nu}_n \rightsquigarrow \bar{\nu} \qquad in \ l_\infty([0,1]^2 \times \mathcal{Z}),$$

*where $\bar{\nu}$ is a centered Gaussian process whose covariance kernel is given by*

$$c(r_1, r_2) \equiv \begin{cases} \langle \beta_{u_1}, \beta_{u_2} \rangle \langle \gamma_{y_1}^\perp, \gamma_{y_2}^\perp \rangle, & \text{if } z_1 = z_2, \\ 0, & \text{if } z_1 \neq z_2. \end{cases}$$

Theorem 3 shows that we can generate distribution-free tests based on $\bar{\nu}_n$. Test statistics are constructed using an appropriate functional $\Gamma$: for example,

$$\Gamma \bar{\nu}_n = \sup_{(u,y,z) \in [0,1]^2 \times \mathcal{Z}} |\bar{\nu}_n(u,y,z)|$$

or

$$\Gamma \bar{\nu}_n = \left( \sum_{z \in \mathcal{Z}} \int_{[0,1]^2} \bar{\nu}_n(u,y,z)^2 \, d(u,y) \right)^{1/2}.$$

When the Gaussian process $\bar{\nu}$ can be simulated, asymptotic critical values can be read from the distribution of $\Gamma \bar{\nu}$. When this is not possible or difficult, one may consider the following bootstrap procedure. Take $\omega_{i,b}$ as in Section 4. Define the bootstrap process

$$\bar{\nu}_{n,b}^*(u,y,z) \equiv \frac{1}{\sqrt{n}} \sum_{i=1}^n \frac{\omega_{i,b} \beta_u(\hat{U}_i) \{ 1\{Z_i = z\} - \hat{p}_{z,i}(\hat{U}_i) \} \gamma_y^\perp(\hat{Y}_i)}{\sqrt{\hat{p}_{z,i}(\hat{U}_i) - \hat{p}_{z,i}(\hat{U}_i)^2}}.$$

We construct the bootstrap test statistics $\bar{T}_{n,b}^* = \Gamma \bar{\nu}_{n,b}^*, b = 1, \ldots, B$, using an appropriate functional $\Gamma$.

THEOREM 4. *Suppose that the conditions of Theorem 3 hold. Then under* $H_0$,

$$d(F_{\bar{T}_n^*}, F_{\Gamma \bar{\nu}}) \to 0 \qquad \text{in } P.$$

## 5. Simulation studies.

5.1. *Conditional independence between continuous variables.* We sampled $X_i$ as i.i.d. Unif$[0,1]$ and $Z_i = aX_i + (1-a)\eta_i$, where $\eta_i \sim$ i.i.d. Unif$[0,1]$ and $a \in \{0.2, 0.5\}$.

We first consider the finite sample size properties of the bootstrap tests. For this purpose, $Y_i$'s were generated as follows:

DGP A1: $Y_i = \Phi((X_i - 0.5)/\sqrt{0.2}) + \varepsilon_i,$

DGP A2: $Y_i = \sin(5X_i) + \varepsilon_i,$

where $\varepsilon_i \sim N(0,1)$ and $\Phi$ denotes the standard normal c.d.f. All the DGPs allow $Y_i$ to depend on $Z_i$, but only through $X_i$, and hence belong to the null hypothesis. The DGPs admit different types of nonlinearity in $X_i$.

TABLE 1
*Rejection probabilities under the null hypothesis of conditional independence
among continuous variables*

| DGP | $a$ | $h$ | Exp. | | | Ind. | | |
|-----|-----|-----|------|------|------|------|------|------|
| | | | 1% | 5% | 10% | 1% | 5% | 10% |
| A1 | 0.2 | $0.25 \times n^{-1/5}$ | 0.0140 | 0.0680 | 0.1285 | 0.0210 | 0.0715 | 0.1340 |
| | | $0.50 \times n^{-1/5}$ | 0.0125 | 0.0540 | 0.1165 | 0.0165 | 0.0615 | 0.1160 |
| | | $1.00 \times n^{-1/5}$ | 0.0120 | 0.0525 | 0.1100 | 0.0140 | 0.0585 | 0.1125 |
| | | $2.00 \times n^{-1/5}$ | 0.0115 | 0.0520 | 0.1090 | 0.0125 | 0.0570 | 0.1025 |
| | 0.5 | $0.25 \times n^{-1/5}$ | 0.0225 | 0.0610 | 0.1180 | 0.0160 | 0.0695 | 0.1280 |
| | | $0.50 \times n^{-1/5}$ | 0.0140 | 0.0580 | 0.1090 | 0.0125 | 0.0540 | 0.1065 |
| | | $1.00 \times n^{-1/5}$ | 0.0105 | 0.0540 | 0.0965 | 0.0120 | 0.0455 | 0.1000 |
| | | $2.00 \times n^{-1/5}$ | 0.0195 | 0.0690 | 0.1315 | 0.0170 | 0.0650 | 0.1315 |
| A2 | 0.2 | $0.25 \times n^{-1/5}$ | 0.0200 | 0.0645 | 0.1295 | 0.0155 | 0.0660 | 0.1215 |
| | | $0.50 \times n^{-1/5}$ | 0.0120 | 0.0555 | 0.1055 | 0.0100 | 0.0560 | 0.1130 |
| | | $1.00 \times n^{-1/5}$ | 0.0115 | 0.0500 | 0.1025 | 0.0100 | 0.0490 | 0.1000 |
| | | $2.00 \times n^{-1/5}$ | 0.0120 | 0.0495 | 0.1095 | 0.0130 | 0.0485 | 0.1040 |
| | 0.5 | $0.25 \times n^{-1/5}$ | 0.0240 | 0.0765 | 0.1385 | 0.0185 | 0.0690 | 0.1295 |
| | | $0.50 \times n^{-1/5}$ | 0.0225 | 0.0640 | 0.1200 | 0.0135 | 0.0700 | 0.1155 |
| | | $1.00 \times n^{-1/5}$ | 0.0135 | 0.0635 | 0.1150 | 0.0165 | 0.0495 | 0.1035 |
| | | $2.00 \times n^{-1/5}$ | 0.0760 | 0.2185 | 0.3285 | 0.0225 | 0.0775 | 0.1465 |

We focus on two types of bootstrap-based tests: one with $\beta_u(U) = 1\{U \le u\}$ (denoted "Ind." in the tables) and the other with $\beta_u(U) = \exp(Uu)$ (denoted "Exp." in the tables). Nonparametric estimations in the Rosenblatt transforms were done using kernel estimation with the kernel $K(u) = (15/16)(1-u^2)^2 1\{|u| \le 1\}$. The bandwidths for $\hat{Y}_i$ and $\hat{Z}_i$ were chosen to be the same, being equal to $h = cn^{-1/5}$ with $c$ ranging in $\{0.25, 0.5, 1, 2\}$. In constructing Kolmogorov–Smirnov tests, we used $10^3$ equal-spaced grid points in $[0,1]^3$. The bootstrap Monte Carlo simulation number and the Monte Carlo simulation number for the whole procedure were set to be 2000. The sample size was equal to 100.

Finite sample sizes are reported in Table 1. The rejection probabilities are overall stable over different choices of bandwidths for all the tests, although they are slightly more sensitive to the bandwidth choices in the case of higher correlation between $X_i$ and $Z_i$ (corresponding to $a = 0.5$).

As for the power properties of the tests, we consider the following four data generating processes:

DGP B1: $Y_i = \Phi((X_i - 0.5)/\sqrt{0.2}) + \Phi((Z_i - 0.5)/\sqrt{0.2}) + \varepsilon_i$,

DGP B2: $Y_i = \Phi((X_i - 0.5)/\sqrt{0.2}) + \sin(5Z_i) + \varepsilon_i$,

TABLE 2
*Rejection probabilities under the alternative hypothesis with nominal level 5%*

| $a$ | $h$ | DGP B1 | | DGP B2 | |
|---|---|---|---|---|---|
| | | Exp. | Ind. | Exp. | Ind. |
| 0.2 | $0.25 \times n^{-1/5}$ | 0.7705 | 0.7560 | 0.9960 | 0.9915 |
| | $0.50 \times n^{-1/5}$ | 0.8015 | 0.7795 | 0.9995 | 0.9960 |
| | $1.00 \times n^{-1/5}$ | 0.8250 | 0.8055 | 1.0000 | 0.9975 |
| | $2.00 \times n^{-1/5}$ | 0.8600 | 0.8500 | 1.0000 | 0.9975 |
| 0.5 | $0.25 \times n^{-1/5}$ | 0.4275 | 0.4140 | 0.9400 | 0.7685 |
| | $0.50 \times n^{-1/5}$ | 0.4505 | 0.4440 | 0.9575 | 0.8035 |
| | $1.00 \times n^{-1/5}$ | 0.4785 | 0.4855 | 0.9725 | 0.8340 |
| | $2.00 \times n^{-1/5}$ | 0.7620 | 0.7640 | 0.9695 | 0.8230 |
| | | DGP B3 | | DGP B4 | |
| 0.2 | $0.25 \times n^{-1/5}$ | 0.0895 | 0.1000 | 0.8945 | 0.5745 |
| | $0.50 \times n^{-1/5}$ | 0.0925 | 0.0945 | 0.9285 | 0.6225 |
| | $1.00 \times n^{-1/5}$ | 0.0830 | 0.0910 | 0.9405 | 0.6570 |
| | $2.00 \times n^{-1/5}$ | 0.0485 | 0.0725 | 0.9470 | 0.6850 |
| 0.5 | $0.25 \times n^{-1/5}$ | 0.0680 | 0.0770 | 0.7395 | 0.4480 |
| | $0.50 \times n^{-1/5}$ | 0.0585 | 0.0680 | 0.7935 | 0.4745 |
| | $1.00 \times n^{-1/5}$ | 0.0495 | 0.0650 | 0.8255 | 0.5000 |
| | $2.00 \times n^{-1/5}$ | 0.1270 | 0.0400 | 0.8555 | 0.4815 |

DGP B3: $Y_i = \sin(5X_i) + \Phi((Z_i - 0.5)/\sqrt{0.2}) + \varepsilon_i$,

DGP B4: $Y_i = \Phi((X_i - 0.5)/\sqrt{0.2}) \times \sin(5Z_i) + \varepsilon_i$.

The results are presented in Table 2. We report the results for the nominal level of 5%. The sample sizes were again 100. For all the cases considered, increasing the correlation between $X_i$ and $Z_i$ (changing $a$ from 0.2 to 0.5) decreases the rejection probabilities. This makes sense because as $X_i$ and $Z_i$ become more dependent, the DGP becomes closer to the null hypothesis.

While the rejection probabilities under DGPs B1 and B2 are reasonably high, the rejection probabilities are much higher in the case of DGP B2 than in the case of DGP B1. In the case of DGP B1, $Y_i$ is monotone both in $X_i$ and $Z_i$, and $Z_i$ is linear in $X_i$. Hence conditional on $X_i$, the presence of the term involving $Z_i$ in the regression model is harder to detect than in the case of DGP B2 where $Y_i$ is not monotone in $Z_i$. The results are similar regardless of whether we use $\beta_u(U) = \exp(Uu)$ or $\beta_u(U) = 1\{U \le u\}$ in constructing test statistics. However, interestingly, when the roles of $X_i$ and $Z_i$ are interchanged as in DGP B3, the tests have very weak power. In simulation studies which are not reported here, we found that the empirical power of the tests was around 75%–95% when the component involving $Z_i$ in DGP B3 was taken to be $\sin(5Z_i)$ or $\cos(5Z_i)$ and $a$ was set to be 0.2. Hence, while the type of nonlinearity between $Y_i$ and $X_i$ plays a crucial role for power properties, the properties also significantly hinge on how $Z_i$ is related to $Y_i$.

Under DGP B4, the rejection probabilities are reasonably high. It is also interesting to observe that under DGP B4, the power properties are significantly better for the choice of $\beta_u(U) = \exp(Uu)$ than for the choice of $\beta_u(U) = 1\{U \leq u\}$. This result illustrates the fact that the choice of $\beta_u(U)$ often plays a significant role in determining the power properties of the test.

5.2. *Conditional independence with binary $Z_i$.* Tests of conditional independence in the case of binary $Z_i$ can be used for program evaluations. For example, suppose $Z_i$ is the binary decision of an individual's treatment which depends on the single index of covariates, $\lambda_{\theta_0}(X_i)$. Then, conditional independence $Y_i \perp\!\!\!\perp Z_i | \lambda_\theta(X_i)$ is a testable implication of the absence of treatment effects under the strong ignorability assumption. In the simulation study, we specified the index as $\lambda_{\theta_0}(X_i) = 0.5 \times (\theta_{00} + \theta_{01}X_{1i} + \theta_{02}X_{2i})$, where $X_i = (X_{1i}, X_{2i})$, $X_{1i} \sim \text{Unif}[0,1] + 0.2$ and $X_{2i} \sim \text{Unif}[0,1] - 0.2$. Here $\theta_{01} = \theta_{02} = 1$ and $\theta_{00} = 0$. The treatment decision $Z_i$ was modeled as

$$Z_i = 1\{\lambda_{\theta_0}(X_i) > \eta_i\}, \qquad \eta_i \sim N(0,1).$$

First, we discuss size properties of the bootstrap tests based on $\beta_u(U) = \exp(Uu)$ and on $\beta_u(U) = 1\{U \leq u\}$. For a specification of the null hypothesis, the variable $Y_i$ was specified as

DGP C: $Y_i = 2\Phi(\lambda_{\theta_0}(X_i)/\sqrt{0.2}) + \varepsilon_i, \qquad \varepsilon_i \sim N(0,1).$

For the construction of the test statistic, we first estimated $\theta_0$ using the MLE to obtain $\hat{\theta}$. Using this estimator $\hat{\theta}$, we constructed $\hat{U}_i$. And then we obtained $\hat{p}_{z,i}(\hat{U}_i)$ and $\hat{Y}_i$ using kernel estimation with the kernel $K(u) = (15/16)(1-u^2)^2 1\{|u| \leq 1\}$ as before. As for taking the Kolmogorov–Smirnov functional, we used $20^2$ equal-spaced grid points in $[0,1]^2$.

The results are presented in Table 3. The number $h_1$ represents the bandwidth for $\hat{p}_{z,i}(\hat{U}_i)$ and $h_2$ for $\hat{Y}_i$. The size properties of the tests are fairly good. The rejection probabilities are mostly close to the nominal level, despite the fact that the test statistics involve a multiple number of nonparametric estimators and an empirical probability integral transform and that the sample size is only 100. The performance of the tests is quite stable over the bandwidth choices and is good regardless of the choice of $\beta_u(U) = \exp(Uu)$ or $\beta_u(U) = 1\{U \leq u\}$.

Let us turn to the power properties of the tests. For this, the following specifications in the alternative hypothesis were used:

DGP D1: $Y_i = 0.5\lambda_{\theta_0}(X_i) + \kappa s(Z_i, X_{1i}, X_{2i}) + \varepsilon_i, \qquad \varepsilon_i \sim N(0,1),$

DGP D2: $Y_i = 2\Phi((\lambda_{\theta_0}(X_i) + \kappa s(Z_i, X_{1i}, X_{2i}))/\sqrt{0.2}) + \varepsilon_i, \qquad \varepsilon_i \sim N(0,1),$

TABLE 3
*Rejection probabilities under the null hypothesis when $Z$ is binary*

| | | Exp. | | | Ind. | | |
|---|---|---|---|---|---|---|---|
| **$h_1$** | **$h_2$** | **1%** | **5%** | **10%** | **1%** | **5%** | **10%** |
| $0.25 \times n^{-1/5}$ | $0.25 \times n^{-1/5}$ | 0.0095 | 0.0430 | 0.0920 | 0.0155 | 0.0495 | 0.1025 |
| | $0.50 \times n^{-1/5}$ | 0.0085 | 0.0355 | 0.0875 | 0.0120 | 0.0500 | 0.1000 |
| | $1.00 \times n^{-1/5}$ | 0.0055 | 0.0320 | 0.0735 | 0.0080 | 0.0460 | 0.0975 |
| | $2.00 \times n^{-1/5}$ | 0.0060 | 0.0345 | 0.0815 | 0.0075 | 0.0430 | 0.0930 |
| $0.50 \times n^{-1/5}$ | $0.25 \times n^{-1/5}$ | 0.0115 | 0.0550 | 0.1090 | 0.0150 | 0.0650 | 0.1150 |
| | $0.50 \times n^{-1/5}$ | 0.0110 | 0.0575 | 0.1130 | 0.0135 | 0.0620 | 0.1240 |
| | $1.00 \times n^{-1/5}$ | 0.0145 | 0.0600 | 0.1050 | 0.0160 | 0.0600 | 0.1080 |
| | $2.00 \times n^{-1/5}$ | 0.0100 | 0.0525 | 0.1050 | 0.0170 | 0.0560 | 0.1150 |
| $1.00 \times n^{-1/5}$ | $0.25 \times n^{-1/5}$ | 0.0110 | 0.0555 | 0.1045 | 0.0135 | 0.0590 | 0.1130 |
| | $0.50 \times n^{-1/5}$ | 0.0130 | 0.0535 | 0.1090 | 0.0145 | 0.0575 | 0.1140 |
| | $1.00 \times n^{-1/5}$ | 0.0135 | 0.0540 | 0.1125 | 0.0140 | 0.0550 | 0.1125 |
| | $2.00 \times n^{-1/5}$ | 0.0140 | 0.0540 | 0.1055 | 0.0140 | 0.0620 | 0.1095 |
| $2.00 \times n^{-1/5}$ | $0.25 \times n^{-1/5}$ | 0.0100 | 0.0455 | 0.0930 | 0.0100 | 0.0435 | 0.0935 |
| | $0.50 \times n^{-1/5}$ | 0.0070 | 0.0440 | 0.0960 | 0.0115 | 0.0530 | 0.0965 |
| | $1.00 \times n^{-1/5}$ | 0.0090 | 0.0475 | 0.0990 | 0.0090 | 0.0475 | 0.1015 |
| | $2.00 \times n^{-1/5}$ | 0.0145 | 0.0525 | 0.1080 | 0.0135 | 0.0525 | 0.1080 |

where $s(Z_i, X_{1i}, X_{2i}) = Z_i\{1 + |X_{1i}| + |X_{2i}|\}$. In the example of program evaluations, the second term, $\kappa s(Z_i, X_{1i}, X_{2i})$, accounts for the path the treatment decision affects the outcome $Y_i$ after conditioning on $\lambda_{\theta_0}(X_i)$. This term involves $Z_i$ and the covariate vector $X_i$ nonlinearly. The number $\kappa$ was chosen from $\{0.5, 1\}$.

The rejection probabilities under the alternative hypothesis are presented in Tables 4 and 5. The rejection probabilities against the alternatives DGP D1 are fairly good. It is interesting to note that the rejection probabilities depend on the choice of bandwidths. The performance is almost the same for the choice of $\beta_u(U) = \exp(Uu)$ or $\beta_u(U) = 1\{U \le u\}$.

The numbers in Table 4 show an interesting result that the bandwidth choice for $\hat{p}_{z,i}(\hat{U}_i)$ is more important for the power property of the test than the bandwidth for $\hat{Y}_i$. When there is more smoothing in the estimation of $\hat{p}_{z,i}(\hat{U}_i)$ within the range of bandwidths considered, the rejection probability improves. However, the rejection probabilities are not as sensitive to the bandwidth choices for $\hat{Y}_i$. Similar observations are made for the case with DGP D2, where $Y_i$ relies nonlinearly on the deviation component $\kappa s(Z_i, X_{1i}, X_{2i})$. In this case, the rejection probabilities are mostly better when $\hat{p}_{z,i}(\hat{U}_i)$ involves more smoothing given the range of the bandwidths. However, the nonlinearity has an overall effect of reducing the rejection prob-abilities. It is also interesting to note that the power properties in this case are different between the choices of $\beta_u(U) = \exp(Uu)$ and $\beta_u(U) = 1\{U \le u\}$. The choice of indicator functions yielded a test with better power properties than the choice of exponential functions in this set-up.

TABLE 4
*Rejection probabilities under the alternative hypothesis (DGP D1) when Z is binary:*
*nominal size = 5%*

| | | $\kappa = 0.5$ | | $\kappa = 1.0$ | |
|---|---|---|---|---|---|
| $h_1$ | $h_2$ | Exp. | Ind. | Exp. | Ind. |
| $0.25 \times n^{-1/5}$ | $0.25 \times n^{-1/5}$ | 0.6075 | 0.6955 | 0.9985 | 0.9985 |
| | $0.50 \times n^{-1/5}$ | 0.6090 | 0.7010 | 1.0000 | 1.0000 |
| | $1.00 \times n^{-1/5}$ | 0.6105 | 0.7020 | 1.0000 | 1.0000 |
| | $2.00 \times n^{-1/5}$ | 0.5950 | 0.6850 | 1.0000 | 1.0000 |
| $0.50 \times n^{-1/5}$ | $0.25 \times n^{-1/5}$ | 0.8850 | 0.9205 | 0.6715 | 0.7605 |
| | $0.50 \times n^{-1/5}$ | 0.9015 | 0.9335 | 0.6675 | 0.7590 |
| | $1.00 \times n^{-1/5}$ | 0.9065 | 0.9385 | 0.6695 | 0.7590 |
| | $2.00 \times n^{-1/5}$ | 0.9055 | 0.9365 | 0.6685 | 0.7585 |
| $1.00 \times n^{-1/5}$ | $0.25 \times n^{-1/5}$ | 0.9545 | 0.9455 | 0.9530 | 0.9815 |
| | $0.50 \times n^{-1/5}$ | 0.9690 | 0.9630 | 0.9535 | 0.9815 |
| | $1.00 \times n^{-1/5}$ | 0.9790 | 0.9740 | 0.9540 | 0.9815 |
| | $2.00 \times n^{-1/5}$ | 0.9850 | 0.9850 | 0.9550 | 0.9810 |
| $2.00 \times n^{-1/5}$ | $0.25 \times n^{-1/5}$ | 0.9465 | 0.9320 | 1.0000 | 1.0000 |
| | $0.50 \times n^{-1/5}$ | 0.9665 | 0.9565 | 1.0000 | 1.0000 |
| | $1.00 \times n^{-1/5}$ | 0.9885 | 0.9775 | 1.0000 | 1.0000 |
| | $2.00 \times n^{-1/5}$ | 0.9980 | 0.9965 | 1.0000 | 1.0000 |

To summarize the findings from the simulation results. First, the size properties of the bootstrap methods based on the distribution-free tests are fairly good. Second, the power properties tend to depend on the choice of $\beta_u(U) = \exp(Uu)$ and $\beta_u(U) = 1\{U \le u\}$ as well as on the bandwidth choices. Third, there are alternatives that the tests have only a trivial power. This finding appears to be consistent with the point made by Janssen (2000) that omnibus tests have nearly trivial asymptotic power against all the local alternatives, except for those contained in a finite-dimensional space.

## APPENDIX: PROOFS

Throughout the proofs, the notation $C$ denotes a positive absolute constant, assuming different values in different contexts. For a class $\mathcal{F}$ of measurable functions, $N(\varepsilon, \mathcal{F}, L_q(P))$ and $N_{[\cdot]}(\varepsilon, \mathcal{F}, L_q(P))$ denote the covering and bracketing numbers of $\mathcal{F}$ with respect to the $L_q(P)$-norm [see van der Vaart and Wellner (1996) for their definitions]. Similarly, we define $N(\varepsilon, \mathcal{F}, \|\cdot\|_\infty)$ and $N_{[\cdot]}(\varepsilon, \mathcal{F}, \|\cdot\|_\infty)$ to be the covering and bracketing numbers with respect to $\|\cdot\|_\infty$.

TABLE 5
*Rejection probabilities under the alternative hypothesis (DGP D2) when $Z$ is binary: nominal size $= 5\%$*

| | | $\kappa = 0.5$ | | $\kappa = 1.0$ | |
|---|---|---|---|---|---|
| $h_1$ | $h_2$ | Exp. | Ind. | Exp. | Ind. |
| $0.25 \times n^{-1/5}$ | $0.25 \times n^{-1/5}$ | 0.2200 | 0.2865 | 0.3040 | 0.4175 |
| | $0.50 \times n^{-1/5}$ | 0.2225 | 0.2895 | 0.3020 | 0.4220 |
| | $1.00 \times n^{-1/5}$ | 0.2120 | 0.2870 | 0.3020 | 0.4260 |
| | $2.00 \times n^{-1/5}$ | 0.2150 | 0.2795 | 0.2980 | 0.4160 |
| $0.50 \times n^{-1/5}$ | $0.25 \times n^{-1/5}$ | 0.3245 | 0.3745 | 0.4350 | 0.5560 |
| | $0.50 \times n^{-1/5}$ | 0.3500 | 0.4075 | 0.4675 | 0.5785 |
| | $1.00 \times n^{-1/5}$ | 0.3405 | 0.4100 | 0.4650 | 0.5975 |
| | $2.00 \times n^{-1/5}$ | 0.3460 | 0.4125 | 0.4615 | 0.5915 |
| $1.00 \times n^{-1/5}$ | $0.25 \times n^{-1/5}$ | 0.3490 | 0.3755 | 0.4680 | 0.5685 |
| | $0.50 \times n^{-1/5}$ | 0.3675 | 0.4135 | 0.5020 | 0.6000 |
| | $1.00 \times n^{-1/5}$ | 0.3820 | 0.4420 | 0.5205 | 0.6320 |
| | $2.00 \times n^{-1/5}$ | 0.3925 | 0.4610 | 0.5275 | 0.6525 |
| $2.00 \times n^{-1/5}$ | $0.25 \times n^{-1/5}$ | 0.3275 | 0.3600 | 0.4750 | 0.5535 |
| | $0.50 \times n^{-1/5}$ | 0.3575 | 0.3960 | 0.5065 | 0.5815 |
| | $1.00 \times n^{-1/5}$ | 0.3840 | 0.4445 | 0.5340 | 0.6360 |
| | $2.00 \times n^{-1/5}$ | 0.4425 | 0.5110 | 0.6040 | 0.6995 |

### A.1. Preliminary results.

LEMMA A1. *Let $\Lambda$ be a class of measurable functions such that for each $\lambda \in \Lambda$, $\lambda(X)$ is continuous with a density function (under $P$) bounded by $M > 0$. Let $\mathcal{T}$ be a class of functions of bounded variation that take values in $[-M, M]$. Then for the class $\mathcal{G} \equiv \{\tau \circ \lambda : (\tau, \lambda) \in \mathcal{T} \times \Lambda\}$, it is satisfied that for any $q \geq 1$,*

$$\log N_{[\cdot]}(C_2 \varepsilon, \mathcal{G}, L_q(P)) \leq \log N(\varepsilon^q, \Lambda, \|\cdot\|_\infty) + C_1/\varepsilon,$$

*where $C_1$ and $C_2$ are positive constants depending only on $q$ and $M$.*

PROOF. Let $F_\lambda$ be the c.d.f. of $\lambda(X)$. For any $\lambda_1, \lambda_2 \in \Lambda$,

$$\sup_x |F_{\lambda_1}(\lambda_1(x)) - F_{\lambda_2}(\lambda_2(x))| \leq M\|\lambda_1 - \lambda_2\|_\infty,$$

because the density of $\lambda(X)$, $\lambda \in \Lambda$, is uniformly bounded. From now on, we identify $\lambda(x)$ with $F_\lambda(\lambda(x))$ without loss of generality so that $\lambda(X)$ is uniformly distributed on $[0,1]$. Since $\mathcal{T} \subset \mathcal{T}_+ - \mathcal{T}_-$ where $\mathcal{T}_+$ and $\mathcal{T}_-$ are collections of uniformly bounded, monotone functions, and we can write $\mathcal{T}_+$ and $\mathcal{T}_-$ as unions of increasing functions and decreasing functions, we lose no generality by assuming that each $\tau \in \mathcal{T}$ is decreasing. Hence by the result of Birman and Solomjak (1967), for any $q \geq 1$

$$\tag{13} \log N_{[\cdot]}(\varepsilon, \mathcal{T}, L_q(P)) \leq \frac{C_2}{\varepsilon}$$

for a constant $C_2 > 0$ that does not depend on $P$.

Choose $\{\lambda_1, \ldots, \lambda_{N_1}\}$ such that any $\lambda \in \Lambda$ is assigned with $\lambda_j$ satisfying $\|\lambda_j - \lambda\|_\infty < \varepsilon^q/2$, and take an integer $M_\varepsilon \in [2\varepsilon^{-q}+1, 2\varepsilon^{-q}+2)$ and a set $\{c_1, \ldots, c_{M_\varepsilon}\}$ such that $c_1 = 0$ and

$$c_{m+1} = c_m + \varepsilon^q/2, \qquad m = 1, \ldots, M_\varepsilon - 1,$$

so that $0 = c_1 \leq c_2 \leq \cdots \leq c_{M_\varepsilon - 1} \leq 1 \leq c_{M_\varepsilon}$. Define

$$\tilde{\lambda}_j(x) = c_m \qquad \text{when } \lambda_j(x) \in [c_m, c_{m+1}), \text{ for some } m \in \{1, 2, \ldots, M_\varepsilon - 1\}.$$

For each $j_1 \in \{1, \ldots, N_1\}$, let $P_{j_1}$ be the distribution of $\tilde{\lambda}_{j_1}(X)$ under $P$. Then choose $\{(\tau_k, \Delta_k)\}_{k=1}^{N_2(j_1)}$ such that any $\tau \in \mathcal{T}$ is assigned with a bracket $(\tau_{j_2}, \Delta_{j_2})$ satisfying $|\tau(\bar{\lambda}) - \tau_{j_2}(\bar{\lambda})| \leq \Delta_{j_2}(\bar{\lambda})$ and $\int \Delta_{j_2}(\bar{\lambda})^q P_{j_1}(d\bar{\lambda}) < \varepsilon^q$.

Now, take any $g \equiv \tau \circ \lambda \in \mathcal{G}$ and let $\lambda_{j_1}$ and $\tau_{j_2}$ be such that $\|\lambda_{j_1} - \lambda\|_\infty < \varepsilon^q/2$ and $|\tau - \tau_{j_2}| \leq \Delta_{j_2}$ with $\int \Delta_{j_2}(\bar{\lambda})^q P_{j_1}(d\bar{\lambda}) < \varepsilon^q$. Fix these $j_1$ and $j_2$ and extend the domain of $\Delta_{j_2}$ to $\mathbf{R}$ by setting $\Delta_{j_2}(\bar{\lambda}) = 0$ for all $\bar{\lambda} \in \mathbf{R} \setminus [0,1]$.

Note that

$$\tag{14}
\begin{aligned}
&|g(x) - (\tau_{j_2} \circ \tilde{\lambda}_{j_1})(x)| \\
&\leq |(\tau \circ \lambda)(x) - (\tau \circ \tilde{\lambda}_{j_1})(x)| + |(\tau \circ \tilde{\lambda}_{j_1})(x) - (\tau_{j_2} \circ \tilde{\lambda}_{j_1})(x)| \\
&\leq |(\tau \circ \lambda)(x) - (\tau \circ \tilde{\lambda}_{j_1})(x)| + (\Delta_{j_2} \circ \tilde{\lambda}_{j_1})(x).
\end{aligned}$$

The range of $\tilde{\lambda}_{j_1}$ is finite and $\|\lambda - \tilde{\lambda}_{j_1}\|_\infty \leq \|\lambda - \lambda_{j_1}\|_\infty + \|\lambda_{j_1} - \tilde{\lambda}_{j_1}\|_\infty \leq \varepsilon^q$. Since $\tau$ is decreasing, $|(\tau \circ \lambda)(x) - (\tau \circ \tilde{\lambda}_{j_1})(x)|$ is bounded by $\tau(\tilde{\lambda}_{j_1}(x) - \varepsilon^q) - \tau(\tilde{\lambda}_{j_1}(x) + \varepsilon^q)$, or by

$$\tau_{j_2}(\tilde{\lambda}_{j_1}(x) - \varepsilon^q) - \tau_{j_2}(\tilde{\lambda}_{j_1}(x) + \varepsilon^q) + \Delta_{j_2}(\tilde{\lambda}_{j_1}(x) - \varepsilon^q) + \Delta_{j_2}(\tilde{\lambda}_{j_1}(x) + \varepsilon^q).$$

Write the difference $\tau_{j_2}(\tilde{\lambda}_{j_1}(x) - \varepsilon^q) - \tau_{j_2}(\tilde{\lambda}_{j_1}(x) + \varepsilon^q)$ as $A_1(x) + A_2(x) + A_3(x) + A_4(x)$, where

$$
\begin{aligned}
A_1(x) &\equiv \tau_{j_2}(\tilde{\lambda}_{j_1}(x) - \varepsilon^q) - \tau_{j_2}(\tilde{\lambda}_{j_1}(x) - \varepsilon^q/2), \\
A_2(x) &\equiv \tau_{j_2}(\tilde{\lambda}_{j_1}(x) - \varepsilon^q/2) - \tau_{j_2}(\tilde{\lambda}_{j_1}(x)), \\
A_3(x) &\equiv \tau_{j_2}(\tilde{\lambda}_{j_1}(x)) - \tau_{j_2}(\tilde{\lambda}_{j_1}(x) + \varepsilon^q/2)
\end{aligned}
$$

and
$$A_4(x) \equiv \tau_{j_2}(\tilde{\lambda}_{j_1}(x) + \varepsilon^q/2) - \tau_{j_2}(\tilde{\lambda}_{j_1}(x) + \varepsilon^q).$$

Due to the construction of $\tilde{\lambda}_{j_1}(x)$, we write $A_1(x)$ as
$$\sum_{m=1}^{M_\varepsilon-1} \{\tau_m^U(j_2) - \tau_m^L(j_2)\} \times 1\{c_m \leq \lambda_{j_1}(x) < c_{m+1}\},$$

where $\tau_m^U(j_2) = \tau_{j_2}(c_m - \varepsilon^q)$ and $\tau_m^L(j_2) = \tau_{j_2}(c_m - \varepsilon^q/2)$. Since $\tau_{j_2}$ is decreasing, $\tau_m^L(j_2) \leq \tau_m^U(j_2)$, and since $c_{m+1} = c_m + \varepsilon^q/2$,
$$\tau_{m+1}^U(j_2) = \tau_{j_2}(c_{m+1} - \varepsilon^q) = \tau_{j_2}(c_m - \varepsilon^q/2) = \tau_m^L(j_2), \qquad m = 1, \ldots, M_\varepsilon - 1.$$

Hence, we conclude
$$\tau_{M_\varepsilon-1}^L(j_2) \leq \cdots \leq \tau_{m+1}^U(j_2) = \tau_m^L(j_2) \leq \tau_m^U(j_2) = \tau_{m-1}^L(j_2) \leq \cdots \leq \tau_1^U(j_2).$$

Suppose that $\tau_1^U(j_2) = \tau_{M_\varepsilon-1}^L(j_2)$. Then $A_1(x) = 0$ and the $L_q(P)$-norm of $A_1$ is trivially zero. Suppose that $\tau_1^U(j_2) > \tau_{M_\varepsilon-1}^L(j_2)$. Since $\tau_{j_2}$ is uniformly bounded, we have $\tau_1^U(j_2) - \tau_{M_\varepsilon-1}^L(j_2) < C < \infty$ for some $C > 0$. Define
$$\tilde{\Delta}_{j_1,j_2}(x) \equiv \sum_{m=1}^{M_\varepsilon-1} \frac{\tau_m^U(j_2) - \tau_m^L(j_2)}{\tau_1^U(j_2) - \tau_{M_\varepsilon-1}^L(j_2)} \times 1\{c_m \leq \lambda_{j_1}(x) < c_{m+1}\}.$$

Let $p_m(j_1) \equiv P\{c_m \leq \lambda_{j_1}(X) < c_{m+1}\}$. Since $\tilde{\Delta}_{j_1,j_2}(x) \leq 1$, $\tilde{\Delta}_{j_1,j_2}^q(x) \leq \tilde{\Delta}_{j_1,j_2}(x)$ so that
$$\mathbf{E}\tilde{\Delta}_{j_1,j_2}^q(X) \leq \sum_{m=1}^{M_\varepsilon-1} \frac{\tau_m^U(j_2) - \tau_m^L(j_2)}{\tau_1^U(j_2) - \tau_{M_\varepsilon-1}^L(j_2)} \times p_m(j_1) \leq \varepsilon^q/2,$$

because $p_m(j_1) \leq \varepsilon^q/2$ for $m \in \{1, \ldots, M_\varepsilon - 1\}$. Thus the $L_q(P)$-norm of $A_1$ is bounded by $C\varepsilon$. We can deal with the functions $A_j(x)$, $j = 2, 3, 4$, similarly by redefining $\tau_m^U(j_2)$ and $\tau_m^L(j_2)$.

From (14) we can bound $|g(x) - (\tau_{j_2} \circ \tilde{\lambda}_{j_1})(x)|$ by

$$\begin{aligned}
(15) \quad & A_1(x) + A_2(x) + A_3(x) + A_4(x) \\
& + (\Delta_{j_2} \circ \tilde{\lambda}_{j_1})(x) + \Delta_{j_2}(\tilde{\lambda}_{j_1}(x) - \varepsilon^q) + \Delta_{j_2}(\tilde{\lambda}_{j_1}(x) + \varepsilon^q) \\
& \equiv \Delta_{j_1,j_2}^*(x).
\end{aligned}$$

Now, let us compute $[\mathbf{E}\{\Delta_{j_1,j_2}^*(X)\}^q]^{1/q}$. The $L_q(P)$-norm of the first four functions is bounded by $C\varepsilon$, as we proved before. By the choice of $\Delta_{j_2}$,
$$\mathbf{E}[\Delta_{j_2}^q(\tilde{\lambda}_{j_1}(X))] = \int \Delta_{j_2}(\bar{\lambda})^q P_{j_1}(d\bar{\lambda}) < \varepsilon^q.$$

Let us turn to the last two terms in (15). Note that

$$\mathbf{E}[\Delta_{j_2}^q(\tilde{\lambda}_{j_1}(X) - \varepsilon^q)] = \sum_{m=1}^{M_\varepsilon-1} \Delta_{j_2}^q(c_m - \varepsilon^q)p_m(j_1)$$

$$\leq \sum_{m=1}^{M_\varepsilon-1} \Delta_{j_2}^q(c_m - \varepsilon^q)\varepsilon^q/2 = \sum_{m=1}^{M_\varepsilon-3} \Delta_{j_2}^q(c_m)\varepsilon^q/2$$

$$\leq \mathbf{E}[\Delta_{j_2}^q(\tilde{\lambda}_{j_1}(X))] < \varepsilon^q.$$

The second equality is due to our setting $\Delta_{j_2}(c) = 0$ for $c \in \mathbf{R} \setminus [0,1]$. The second to the last inequality is due to the fact that $p_m(j_1) = \varepsilon^q/2$ for all $m \in \{1, \ldots, M_\varepsilon - 2\}$. Similarly, $\mathbf{E}[\Delta_{j_2}^q(\tilde{\lambda}_{j_1}(X) + \varepsilon^q)] < \varepsilon^q$. Combining these results, $\mathbf{E}[\{\Delta_{j_1,j_2}^*(X)\}^q] \leq C_1^q \varepsilon^q$, for some constant $C_1 > 0$, yielding the result that

$$\log N_{[\cdot]}(C_1\varepsilon, \mathcal{G}, L_q(P)) \leq \log N(\varepsilon^q/2, \Lambda, \|\cdot\|_\infty) + C_2/\varepsilon.$$

With redefinitions of constants and $\varepsilon$, this completes the proof. $\quad\square$

Given $x_{(n)} \equiv (x_j)_{j=1}^n \in \mathbf{R}^{nd_X}$, let

$$(16) \qquad G_{n,\lambda}(\cdot; x_{(n)}) \equiv \frac{1}{n}\sum_{j=1}^n \mathbf{1}\{\lambda(x_j) \leq \cdot\}.$$

Lemma A1 yields the following bracketing entropy bound by taking $\tau = \beta_u \circ G_{n,\lambda}$. Certainly, this $\tau$ is bounded variation, because $G_{n,\lambda}$ is increasing.

COROLLARY A1. *Let $\Lambda$ and $M$ be as in Lemma A1, and for $\beta_u$ in Assumption 2(i) let*

$$\mathcal{B}_n \equiv \{\beta_u(G_{n,\lambda}(\lambda(\cdot); x_{(n)})) : (u, \lambda, x_{(n)}) \in [0,1] \times \Lambda \times \mathbf{R}^{nd_X}\}.$$

*Then $\log N_{[\cdot]}(C_2\varepsilon, \mathcal{B}_n, L_q(P)) \leq \log N(\varepsilon^q, \Lambda, \|\cdot\|_\infty) + C_1/\varepsilon$, for any $q \geq 1$, where $C_1$ and $C_2$ are positive constants depending only on $q$ and $M$.*

The following lemma is useful for establishing a bracketing entropy bound of a class in which conditional c.d.f. estimators realize.

LEMMA A2. *We introduce three classes of functions. First, let $\mathcal{F}_n$ be a sequence of classes of maps $\phi(\cdot, \cdot) : \mathbf{R} \times \mathcal{S}_n \to [0,1]$, such that (a) $\mathcal{S}_n \subset [-s_n, s_n]$, $s_n > 1$, (b) for each $v \in \mathcal{S}_n$, $\phi(\cdot, v)$ is bounded variation and (c) for each $\varepsilon > 0$,*

$$(17) \qquad \sup_{(y,v) \in \mathbf{R} \times \mathcal{S}_n} \sup_{\eta \in [0,\varepsilon]} |\phi(y, v + \eta) - \phi(y, v - \eta)| < M_n\varepsilon$$

for some sequence $M_n > 1$. Second, let $\mathcal{G}$ be a class of measurable functions $G : \mathbf{R}^{d_X} \to \mathcal{S}_n$. Finally, let $\mathcal{J}_n^{\mathcal{G}} \equiv \{ \phi(\cdot, G(\cdot)) : (\phi, G) \in \mathcal{F}_n \times \mathcal{G} \}$.

Then for any probability $P$, and for any $q > 1$,

$$\log N_{[\cdot]}(\varepsilon, \mathcal{J}_n^{\mathcal{G}}, L_q(P))$$
$$\leq C + C \log(2 s_n M_n + \varepsilon)$$
$$+ C \{ (2 s_n M_n + \varepsilon)^{1/q} \varepsilon^{-(q+1)/q} - \log(\varepsilon) + \log N_{[\cdot]}(C \varepsilon / M_n, \mathcal{G}, L_q(P)) \}$$

for some $C > 0$.

PROOF. Fix $\varepsilon > 0$ and take a partition $\mathcal{S}_n = \bigcup_{k=1}^{J_\varepsilon} B(b_k)$ where $B(b_k)$ is a set contained in an $\varepsilon / M_n$-interval centered at $b_k$, and $J_\varepsilon \leq 2 s_n M_n / \varepsilon + 1$. Let $\mathcal{F}_n(b_k) \equiv \{ \phi(\cdot, b_k) : \phi \in \mathcal{F}_n \}$. For each $b_k$, take $\{ (f_{k,j}, \Delta_{k,j}) \}_{j=1}^{N_\varepsilon(b_k)}$ such that to any $f \in \mathcal{F}_n(b_k)$ corresponds $(f_{k,j}, \Delta_{k,j})$ such that $|f(y) - f_{k,j}(y)| \leq \Delta_{k,j}(y)$ and $\int \Delta_{k,j}^q(y) P(dy) \leq \varepsilon^{q+1} / (2 s_n M_n + \varepsilon)$.

Given $\phi \in \mathcal{F}_n$, we let $\tilde{f}_j(y, v) \equiv \sum_{k=1}^{J_\varepsilon} f_{k,j}(y) A_k(v)$, where $A_k(v) \equiv \mathbf{1} \{ v \in B(b_k) \}$ and $f_{k,j}(y)$ is such that $|\phi(y, b_k) - f_{k,j}(y)| \leq \Delta_{k,j}(y)$ and

$$\int \Delta_{k,j}^q(y) P(dy) \leq \varepsilon^{q+1} / (2 s_n M_n + \varepsilon).$$

Since $\mathcal{F}_n(b_k)$ is a uniformly-bounded class of bounded variation functions, the smallest number of such $(k, j)$'s are bounded by

$$(18) \quad J_\varepsilon \exp \left( \frac{C(2 s_n M_n + \varepsilon)^{1/q}}{\varepsilon^{(q+1)/q}} \right) \leq \left( \frac{2 s_n M_n}{\varepsilon} + 1 \right) \exp \left( \frac{C(2 s_n M_n + \varepsilon)^{1/q}}{\varepsilon^{(q+1)/q}} \right).$$

Then we bound $|\phi(y, v) - \tilde{f}_j(y, v)|$ by

$$(19) \quad \left| \sum_{k=1}^{J_\varepsilon} \{ \phi(y, v) - \phi(y, b_k) \} A_k(v) \right| + \left| \sum_{k=1}^{J_\varepsilon} \{ \phi(y, b_k) - f_{k,j}(y) \} A_k(v) \right|.$$

The second term is bounded by $\sum_{k=1}^{J_\varepsilon} \Delta_{k,j}(y) A_k(v)$, and the first term, by $\varepsilon$, due to (17). Hence

$$|\phi(y, v) - \tilde{f}_j(y, v)| \leq \bar{\Delta}_j(y, v),$$

where $\bar{\Delta}_j(y, v) \equiv \sum_{k=1}^{J_\varepsilon} \Delta_{k,j}(y) A_k(v) + \varepsilon$. By the Hölder inequality,

$$\left( \sum_{k=1}^{J_\varepsilon} \Delta_{k,j}(y) A_k(v) \right)^q \leq \sum_{k=1}^{J_\varepsilon} \Delta_{k,j}^q(y).$$

Hence, $\int \bar{\Delta}_j^q(y, v) P(dy, dv)$ is bounded by

$$C \sum_{k=1}^{J_\varepsilon} \int \Delta_{k,j}^q(y) P(dy) + C \varepsilon^q \leq \frac{C J_\varepsilon \varepsilon^{q+1}}{2 s_n M_n + \varepsilon} + C \varepsilon^q \leq 2 C \varepsilon^q,$$

yielding the inequality [from (18)]

$$
\begin{aligned}
(20) \quad & \log N_{[\cdot]}(C\varepsilon, \mathcal{F}_n, L_q(P)) \\
& \leq C + \log(2s_n M_n + \varepsilon) + C(2s_n M_n + \varepsilon)^{1/q} / \varepsilon^{(q+1)/q} - \log(\varepsilon).
\end{aligned}
$$

Now, take $(G_k, \Delta_k)_{k=1}^{N_1}$ such that to any $G \in \mathcal{G}$ corresponds $(G_j, \Delta_j)$ such that $|G - G_j| \leq \Delta_j$ and $\mathbf{E}\Delta_j^q(X) < (\varepsilon/M_n)^q$. Take $(\phi_k, \tilde{\Delta}_k)_{k=1}^{N_2(j)}$ such that for any $\phi \in \mathcal{F}_n$, there exists $(\phi_k, \tilde{\Delta}_k)$ such that $|\phi(y, v) - \phi_k(y, v)| < \tilde{\Delta}_k(y, v)$ and $\mathbf{E}\tilde{\Delta}_k^q(Y, G_j(X)) < \varepsilon^q$. Now, we bound $|\phi(y, G(x)) - \phi_k(y, G_j(x))|$ by

$$
\begin{aligned}
& |\phi(y, G(x)) - \phi(y, G_j(x))| + |\phi(y, G_j(x)) - \phi_k(y, G_j(x))| \\
& \qquad \leq M_n \Delta_j(x) + \tilde{\Delta}_k(y, G_j(x)) \equiv \bar{\Delta}_{j,k}(y, x).
\end{aligned}
$$

Since $\mathbf{E}\bar{\Delta}_{j,k}^q(Y, X) \leq C\varepsilon^q$, we conclude that

$$
\begin{aligned}
& \log N_{[\cdot]}(C\varepsilon, \mathcal{J}_n^{\mathcal{G}}, L_q(P)) \\
& \qquad \leq C \log N_{[\cdot]}(\varepsilon, \mathcal{F}_n, L_q(P)) + C \log N_{[\cdot]}(\varepsilon/M_n, \mathcal{G}, L_q(P)).
\end{aligned}
$$

Combined with (20), this yields the desired result. $\square$

Fix $\lambda_0 : \mathbf{R}^{d_X} \to \mathbf{R}$ in a uniformly bounded class $\Lambda$, and let $F_0$ be the distribution function of $\lambda_0(X)$. Let $F_{n,0,i}$ and $F_{n,\lambda,i}$ be the empirical distribution functions of $\{\lambda_0(X_j)\}_{j=1, j\neq i}^n$ and $\{\lambda(X_j)\}_{j=1, j\neq i}^n$, $\{X_j\}_{j=1}^n$ being i.i.d. $\sim P$.

LEMMA A3. *Define $\Lambda_n \equiv \{\lambda \in \Lambda : \|\lambda - \lambda_0\|_\infty \leq Cn^{-1/2}\}$ for $C > 0$ and assume that $\Lambda_n$ satisfy the conditions for $\Lambda$ in Lemma A1 and that $\log N(C\varepsilon, \Lambda_n, \|\cdot\|_\infty) < C\varepsilon^{-r_\Lambda}$ for some $C > 0$ and $r_\Lambda \in [0, 1)$. Then,*

$$
\mathbf{E}\Big[\max_{1 \leq i \leq n} \sup_{\lambda \in \Lambda_n} \sup_{x \in \mathbf{R}^{d_X}} |F_{n,\lambda,i}(\lambda(x)) - F_0(\lambda_0(x))|\Big] = O(n^{-1/2}).
$$

PROOF. Let every $\lambda$ in $\Lambda$ be bounded in $[-M, M]$. It suffices to show that

$$
\begin{aligned}
(21) \quad & \sup_{\lambda \in \Lambda_n} \sup_{x \in \mathbf{R}^{d_X}} |F_\lambda(\lambda(x)) - F_0(\lambda_0(x))| = O(n^{-1/2}) \quad \text{and} \\
& \mathbf{E}\Big[\max_{1 \leq i \leq n} \sup_{\lambda \in \Lambda_n} \sup_{\bar{\lambda} \in [-M, M]} |F_{n,\lambda,i}(\bar{\lambda}) - F_\lambda(\bar{\lambda})|\Big] = O(n^{-1/2}).
\end{aligned}
$$

The first statement of (21) follows, due to $\lambda_0(X)$ having a bounded density and $\sup_{\lambda \in \Lambda_n} \|\lambda - \lambda_0\|_\infty = O(n^{-1/2})$ by the definition of $\Lambda_n$. As for the second statement, write

$$
|F_{n,\lambda,i}(\bar{\lambda}) - F_\lambda(\bar{\lambda})| \leq \left| \frac{1}{n-1} \sum_{j=1, j\neq i}^n (1\{\lambda(X_j) \leq \bar{\lambda}\} - P\{\lambda(X_j) \leq \bar{\lambda}\}) \right|.
$$

For any $\lambda_1, \lambda_2 \in \Lambda_n$ and $\bar\lambda_1, \bar\lambda_2 \in [-M, M]$,

$$\mathbf{E}|1\{\lambda_1(X_j) \le \bar\lambda_1\} - 1\{\lambda_2(X_j) \le \bar\lambda_2\}| \le C\{\|\lambda_1 - \lambda_2\|_\infty + |\bar\lambda_1 - \bar\lambda_2|\}.$$

This implies that

$$\log N_{[\cdot]}(\varepsilon, \mathcal{F}_n, L_2(P)) \le \log N(C\varepsilon^2, \Lambda_n \times [-M, M], \|\cdot\|_\infty + |\cdot|) \le C\varepsilon^{-2r_\Lambda},$$

where $\mathcal{F}_n \equiv \{1\{\lambda(\cdot) \le \bar\lambda\} : (\lambda, \bar\lambda) \in \Lambda \times [-M, M]\}$. Since $r_\Lambda < 1$, the second result of (21) follows from the maximal inequality of Pollard (1989) [see Theorem A.2 of van der Vaart (1996)]. $\square$

Let $\Lambda_n, \lambda_0, F_0$, and $F_{n,\lambda,i}$ be as in Lemma A3 and define

$$U_{n,\lambda,i} \equiv F_{n,\lambda,i}(\lambda(X_i)) \quad \text{and} \quad U \equiv F_0(\lambda_0(X))$$

for $\lambda \in \Lambda_n$. Let $\mathcal{S}_W$ be a subset of $\mathbf{R}^{d_W}$ and introduce $\Phi_n$, a class of functions $\varphi : \mathcal{S}_W \to \mathbf{R}$. Then we define a kernel estimator of $g_\varphi(u) \equiv \mathbf{E}[\varphi(W)|U = u]$:

$$(22) \qquad \hat{g}_{\varphi,\lambda,i}(u) \equiv \frac{\sum_{j=1, j \ne i}^n \varphi(W_j) K_h(U_{n,\lambda,j} - u)}{\sum_{j=1, j \ne i}^n K_h(U_{n,\lambda,j} - u)}.$$

The following lemma proves uniform convergence of $\hat{g}_{\varphi,\lambda,i}(u)$ over $(\varphi, \lambda) \in \Phi_n \times \Lambda_n$. A related result without the transform $F_{n,\lambda,i}$ was obtained by Andrews (1995).

LEMMA A4.   *Let $\Lambda_n$ be as in Lemma A3. Furthermore, assume the following:*

(A1) $\log N(C\varepsilon, \Lambda_n, \|\cdot\|_\infty) < C\varepsilon^{-r_\Lambda}$ *for some $C > 0$ and $r_\Lambda \in [0, 1)$.*

(A2) $\Phi_n$ *has an envelope $\bar\varphi$ such that $\|\bar\varphi\|_{P,2} < \infty$,*

$$\sup_{u \in [0,1]} \mathbf{E}[|\bar\varphi(W_i)| | U_i = u] < \infty,$$

*and for some $r_\Phi \in [0, 2)$, $\log N_{[\cdot]}(\varepsilon, \Phi_n, L_2(P)) < \varepsilon^{-r_\Phi}$.*

(A3) $g_\varphi(\cdot), \varphi \in \Phi_n$, *is twice continuously differentiable with uniformly bounded derivatives.*

(A4) $K$ *is symmetric, nonnegative, continuously differentiable, has a compact support, and $\int_{-\infty}^\infty K(s) = 1$.*

*Then as $h \to 0$ and $n^{-1/2}h^{-1} \to 0$ with $n \to \infty$,*

$$\max_{1 \le i \le n} \sup_{(\varphi, \lambda) \in \Phi_n \times \Lambda_n} |\hat{g}_{\varphi,\lambda,i}(u) - g_\varphi(u)|$$

$$\le \Delta_n(u) + O(n^{-1/2}h^{-1}\sqrt{-\log h}),$$

*where $\Delta_n(u) = C(h^2 1\{|u - 1| > h\} + h 1\{|u - 1| \le h\})$ for some $C > 0$.*

*Furthermore, if $\Phi_n$ is uniformly bounded,*

$$\max_{1 \le i \le n} \sup_{(\varphi, \lambda) \in \Phi_n \times \Lambda_n} \sup_{u \in [0,1]} |\hat{g}_{\varphi,\lambda,i}(u) - g_\varphi(u)|$$

$$\le \Delta_n(u) + O(n^{-1/2}h^{-1}).$$

PROOF. We consider the first statement. Without loss of generality, assume that $K$ has a support in $[-1, 1]$, and $K(\cdot) \leq M$ for some $M > 0$. By Lemma A3, it suffices to focus on the event $\sup_x |F_{n,\lambda,i}(\lambda(x)) - F_0(\lambda_0(x))| < Mn^{-1/2}$ for large $M > 0$. Define

$$\hat{\rho}_{\varphi,\lambda,i}(u) \equiv \frac{1}{n-1} \sum_{j=1, j \neq i}^{n} K_h(U_{n,\lambda,j} - u)\varphi(W_j)$$

and

$$\hat{f}_{\lambda,i}(u) \equiv \frac{1}{n-1} \sum_{j=1, j \neq i}^{n} K_h(U_{n,\lambda,j} - u).$$

First, let us prove that uniformly over $1 \leq i \leq n$,

$$
\begin{aligned}
(23) \quad & \sup_{(\varphi,\lambda) \in \Phi_n \times \Lambda_n} |\hat{\rho}_{\varphi,\lambda,i}(u) - \xi_{1n}(u)g_\varphi(u)| \\
& \leq \Delta_n(u) + O(n^{-1/2}h^{-1}\sqrt{-\log h}),
\end{aligned}
$$

where $\xi_{1n}(u) = \int_{(-u/h)\vee(-1)}^{((1-u)/h)\wedge 1} K(v)\,dv$. Write $\hat{\rho}_{\varphi,\lambda,i}(u) - \xi_{1n}(u)g_\varphi(u)$ as

$$
\begin{aligned}
& \frac{1}{n-1} \sum_{j=1, j \neq i}^{n} \{K_h(U_{n,\lambda,j} - u) - K_h(U_j - u)\}\varphi(W_j) \\
& \quad + \frac{1}{n-1} \sum_{j=1, j \neq i}^{n} K_h(U_j - u)\varphi(W_j) - \xi_{1n}(u)g_\varphi(u) \\
& \equiv A_{1n} + A_{2n}.
\end{aligned}
$$

We deal with $A_{1n}$ following the proof of Lemma 4.5 of Stute and Zhu (2005). Since $K$ has support in $[-1, 1]$, $A_{1n}$ is the sum of those $j$'s such that either $|U_{n,\lambda,j} - u| \leq h$ or $|U_j - u| \leq h$. By Lemma A3,

$$\max_{1 \leq j \leq n, \lambda \in \Lambda_n} \sup |U_{n,\lambda,j} - U_j| = O_P(n^{-1/2}).$$

Since $n^{-1/2}h^{-1} \to 0$, $A_{1n}$ is equal to

$$(24) \quad \frac{1}{(n-1)h^2} \sum_{j=1, j \neq i}^{n} K'(\Delta_{ij})\varphi(W_j)\{U_{n,\lambda,j} - U_j\}1\{|U_j - u| \leq Ch\}$$

with probability approaching one, where $\Delta_{ij}$ lies between $(U_{n,\lambda,j} - u)/h$ and $(U_j - u)/h$. Therefore, by (A4), the term in (24) is bounded by

$$O_P(n^{-1/2}) \times \frac{1}{(n-1)h^2} \sum_{j=1, j \neq i}^{n} |\bar{\varphi}(W_j)|1\{|U_j - u| \leq Ch\}$$

with probability approaching one. By (A2), the expectation of the above sum is $O(h^{-1})$, yielding

$$(25) \qquad A_{1n} = O_P(n^{-1/2}h^{-1}).$$

This leaves us to deal with $A_{2n}$ which we write as

$$A_{2n} = \frac{1}{n-1} \sum_{j=1, j \neq i}^{n} K_h(U_j - u)\{\varphi(W_j) - g_\varphi(U_j)\}$$

$$+ \frac{1}{n-1} \sum_{j=1, j \neq i}^{n} K_h(U_j - u)g_\varphi(U_j) - \xi_{1n}(u)g_\varphi(u)$$

$$\equiv \text{(I)} + \text{(II)}.$$

The sum (I) is a mean-zero process. Define $\phi_u(y, v) \equiv K((y - u)/h)v$, $\tilde{\varphi}(w, u; \varphi) \equiv \varphi(w) - g_\varphi(u)$, and

$$\tilde{\Phi}_n \equiv \{\tilde{\varphi}(\cdot, \cdot; \varphi) : \varphi \in \Phi_n\}$$

and

$$\mathcal{J}_{1n} \equiv \{\phi_u(\cdot, \tilde{\varphi}(\cdot)) : (u, \tilde{\varphi}) \in [0, 1] \times \tilde{\Phi}_n\}.$$

Take $J(w, u) = M\{\tilde{\varphi}(w) + g_{\tilde{\varphi}}(u)\}$ as an envelop for $\mathcal{J}_{1n}$. Then, $\mathbf{E}|J(W, U)|^2 < \infty$. Note that

$$|\phi_{u_1}(y, \tilde{\varphi}(w, u; \varphi)) - \phi_{u_2}(y, \tilde{\varphi}(w, u; \varphi))| \leq C|u_1 - u_2||\tilde{\varphi}(w) + g_{\tilde{\varphi}}(u)|/h$$

and

$$|\phi_u(y, \tilde{\varphi}(w, u; \varphi_1)) - \phi_u(y, \tilde{\varphi}(w, u; \varphi_2))|$$
$$\leq C|\varphi_1(w) - \varphi_2(w)| + C|g_{\varphi_1}(u) - g_{\varphi_2}(u)|.$$

Since $\mathbf{E}|g_{\varphi_1}(U) - g_{\varphi_2}(U)|^q \leq \mathbf{E}|\varphi_1(W) - \varphi_2(W)|^q$, $q \geq 1$, we conclude that

$$\log N_{[\cdot]}(\varepsilon, \mathcal{J}_{1n}, L_2(P))$$
$$(26) \qquad \leq \log N(C\varepsilon h, [0, 1], |\cdot|) + \log N(C\varepsilon, \Phi_n, L_2(P))$$
$$\leq C - (\log \varepsilon + \log h) + C\varepsilon^{-r_\Phi}.$$

Therefore, by the maximal inequality of Pollard ([1989]) [e.g., Theorem A.2 of van der Vaart ([1996])],

$$\mathbf{E}\left[ \sup_{(u,\varphi) \in [0,1] \times \Phi_n} \left| \frac{1}{n-1} \sum_{j=1, j \neq i}^{n} K\left( \frac{U_j - u}{h} \right)\{\varphi(W_j) - g_\varphi(U_j)\} \right| \right]$$

$$(27) \qquad \leq \frac{C}{\sqrt{n-1}} \int_0^C \sqrt{1 + \log N_{[\cdot]}(\varepsilon, \mathcal{J}_{1n}, L_2(P))}\, d\varepsilon$$

$$= O(n^{-1/2}\sqrt{-\log h}).$$

Hence the uniform convergence rate for (I) is $O(n^{-1/2}h^{-1}\sqrt{-\log h})$.

Let us turn to (II). For this, write (II) as

$$
\begin{aligned}
(28) \quad & \frac{1}{n-1}\sum_{j=1,j\neq i}^{n}\{K_h(U_j-u)g_\varphi(U_j)-\mathbf{E}[K_h(U_j-u)g_\varphi(U_j)]\} \\
& + \mathbf{E}[K_h(U_j-u)g_\varphi(U_j)] - \xi_{1n}(u)g_\varphi(u).
\end{aligned}
$$

From the steps to prove (27), the first sum is $O_P(n^{-1/2}h^{-1}\sqrt{-\log h})$ uniformly over $u\in[0,1]$. Write the last difference in (28), for $u\in[0,1]$, as

$$
\frac{1}{h}\int_0^1 K\left(\frac{v-u}{h}\right)g_\varphi(v)\,dv - g_\varphi(u)\xi_{1n}(u) = g_\varphi'(u)\xi_{2n}(u)h + O(h^2),
$$

where $\xi_{2n}(u) = \int_{(-u/h)\vee(-1)}^{((1-u)/h)\wedge 1} vK(v)\,dv$. Therefore,

$$
\text{(II)} \leq \Delta_n(u) + O(n^{-1/2}h^{-1}\sqrt{-\log h}).
$$

Combining (I) and (II), we obtain (23).

Following the proof of (23) with $\varphi=1$, we also obtain

$$
(29) \quad \sup_{(\varphi,\lambda,u)\in\Phi_n\times\Lambda_n\times[0,1]}|\hat{f}_{\lambda,i}(u)-\xi_{1n}(u)| = O(n^{-1/2}h^{-1}\sqrt{-\log h}).
$$

The quantity $|\xi_{1n}(u)-1\{u\in[0,1]\}|$ is bounded by

$$
\max\left\{\int_0^1 K(v)\,dv, \int_{-1}^0 K(v)\,dv\right\} \leq 1/2.
$$

Thus $\xi_{1n}(u)\geq 1/2$ for $u\in[0,1]$. Write $\hat{g}_{\varphi,\lambda,i}(u)-g_\varphi(u)$ as

$$
\hat{\rho}_{\varphi,\lambda,i}(u)\left\{\frac{1}{\hat{f}_{\lambda,i}(u)}-\frac{1}{\xi_{1n}(u)}\right\} + \frac{1}{\xi_{1n}(u)}\{\hat{\rho}_{\varphi,\lambda,i}(u)-\xi_{1n}(u)g_\varphi(u)\}.
$$

We apply (23) and (29) to the first and second terms to obtain the first result of the theorem.

As for the second result, we modify the treatment of (I) above. Note that

$$
\sup_{(y,v)\in\mathbf{R}\times[-s_n,s_n]}\sup_{\eta\in[0,\varepsilon]}|\phi_u(y,v+\eta)-\phi_u(y,v-\eta)| \leq C\varepsilon.
$$

Since $K(\cdot)$ is Lipschitz, it is bounded variation, and so are $K((\cdot-u)/h)$. Since $\Phi_n$ is uniformly bounded, we apply Lemma A2 (with $s_n$ equal to some large $M>0$),

$$
\log N_{[\cdot]}(\varepsilon,\mathcal{J}_{1n},L_2(P)) \leq C + C\{\varepsilon^{-3/2}-\log(\varepsilon)+\varepsilon^{-r_\Phi}\}.
$$

Substituting this bound for the one in (26) and following the same arguments there, we obtain the wanted result. $\quad\square$

Let $\Lambda$ be a uniformly bounded class of real functions and fix $\lambda_0 \in \Lambda$. Let $\mathcal{S}_X$ be the support of a random vector $X$, and $[0,1]^2$ be the support of a random vector $(\tilde{Z}, \tilde{Y})$. Let $f_\lambda(z|y, \bar{\lambda}_1, \bar{\lambda}_0)$ be the conditional density of $\tilde{Z}$ given $\tilde{Y} = y, \lambda_1(X) = \bar{\lambda}_1$ and $\lambda_0(X) = \bar{\lambda}_0$ with respect to a $\sigma$-finite conditional measure. Introduce

$$\mu_\lambda(z|y, \bar{\lambda}_1, \bar{\lambda}_0) \equiv \mathbf{E}[\gamma_z(\tilde{Z})|\tilde{Y} = y, \lambda(X) = \bar{\lambda}_1, \lambda_0(X) = \bar{\lambda}_0]$$

and

$$\mu_\lambda(z|y, \bar{\lambda}_0) \equiv \mathbf{E}[\gamma_z(\tilde{Z})|\tilde{Y} = y, \lambda_0(X) = \bar{\lambda}_0].$$

For each $\lambda \in \Lambda, \mathcal{S}(\lambda)$ be the support of $\lambda(X)$.

LEMMA A5.  *Suppose that there exists $\bar{C} > 0$ such that for each $\lambda_1 \in \Lambda$, each $\delta > 0$ and each $(\bar{\lambda}_1, \bar{\lambda}_0) \in \mathcal{S}(\lambda_1) \times \mathcal{S}(\lambda_0)$,*

$$\sup_{(z, \bar{\lambda}_1) \in [0,1]^2 \times \mathcal{S}(\lambda_1) \,:\, |\bar{\lambda}_1 - \bar{\lambda}_1| < \delta} |f_{\lambda_1}(z|y, \bar{\lambda}_1, \bar{\lambda}_0) - f_{\lambda_1}(z|y, \bar{\lambda}_1, \bar{\lambda}_0)| \le \bar{C}\delta.$$

*Then for each $\delta > 0$ and each $\lambda_1$ in $\Lambda$ such that $\|\lambda_1 - \lambda_0\|_\infty \le \delta$,*

$$\sup_{(z,y,x) \in [0,1]^2 \times \mathcal{S}_X} |\mu_{\lambda_1}(z|y, \lambda_1(x), \lambda_0(x)) - \mu_{\lambda_0}(z|y, \lambda_0(x))| \le 6\bar{C}\delta.$$

PROOF.  Choose $(z, y, x) \in [0,1]^2 \times \mathcal{S}_X$ and $\lambda_1 \in \Lambda$ with $\|\lambda_1 - \lambda_0\|_\infty < \delta$ and let $\bar{\lambda}_0 \equiv \lambda_0(x)$ and $\bar{\lambda}_1 \equiv \lambda_1(x)$. Let $P_y$ be the conditional measure of $(\tilde{Z}, X)$ given $(\tilde{Y}, \lambda_0(X)) = (y, \bar{\lambda}_0)$ and $\mathbf{E}_y$ be the conditional expectation under $P_y$. Let $A_j \equiv 1\{|\lambda_j(X) - \bar{\lambda}_j| \le 3\delta\}, j = 0, 1$. Note that $\mathbf{E}_y[A_0] = 1$ and

$$1 \ge \mathbf{E}_y[A_1] = P\{|\lambda_1(X) - \bar{\lambda}_1| \le 3\delta|Y = y, \lambda_0(X) = \bar{\lambda}_0\}$$
$$\ge P\{|\lambda_0(X) - \bar{\lambda}_0| \le \delta|Y = y, \lambda_0(X) = \bar{\lambda}_0\} = 1.$$

Let $\tilde{\mu}_{\lambda_j}(z|y, \bar{\lambda}_j, \bar{\lambda}_0) \equiv \mathbf{E}_y[\gamma_z(\tilde{Z})A_j]/\mathbf{E}_y[A_j] = \mathbf{E}_y[\gamma_z(\tilde{Z})A_j], j = 0, 1$. Then, for example,

$$|\mu_{\lambda_1}(z|y, \bar{\lambda}_1, \bar{\lambda}_0) - \mu_{\lambda_0}(z|y, \bar{\lambda}_0)|$$
$$\le |\mu_{\lambda_1}(z|y, \bar{\lambda}_1, \bar{\lambda}_0) - \tilde{\mu}_{\lambda_1}(z|y, \bar{\lambda}_1, \bar{\lambda}_0)| + |\tilde{\mu}_{\lambda_1}(z|y, \bar{\lambda}_1, \bar{\lambda}_0) - \mu_{\lambda_0}(z|y, \bar{\lambda}_0)|$$
$$\equiv (\mathrm{I}) + (\mathrm{II}).$$

Let us turn to (I). By the definition of conditional expectation,

$$\tilde{\mu}_{\lambda_1}(z|y, \bar{\lambda}_1, \bar{\lambda}_0) = \int_{\bar{\lambda}_1 - 3\delta}^{\bar{\lambda}_1 + 3\delta} \mu_{\lambda_1}(z|y, \bar{\lambda}, \bar{\lambda}_0)\, dF_{\lambda_1}(\bar{\lambda}|y, \bar{\lambda}_0),$$

where $F_{\lambda_1}(\cdot|y, \bar{\lambda}_0)$ is the conditional c.d.f. of $\lambda_1(X)$ given $(\tilde{Y}, \lambda_0(X)) = (y, \bar{\lambda}_0)$. Because $\int_{\bar{\lambda}_1 - 3\delta}^{\bar{\lambda}_1 + 3\delta} dF_{\lambda_1}(\bar{\lambda}|y, \bar{\lambda}_0) = \mathbf{E}_y[A_1] = 1$ and $|\gamma_z(\tilde{Z})| \le 1$, wp 1,

$$(\mathrm{I}) \le \sup_{v \in [-3\delta, 3\delta] \,:\, \bar{\lambda}_1 + v \in \mathcal{S}(\lambda_1)} |f_{\lambda_1}(z|y, \bar{\lambda}_1 + v, \bar{\lambda}_0) - f_{\lambda_1}(z|y, \bar{\lambda}_1, \bar{\lambda}_0)|.$$

Therefore, by the condition of the lemma, (I) $\leq 3\bar{C}\delta$.

Let us turn to (II) which we write as

$$|\mathbf{E}_y[\gamma_z(\tilde{Z})A_1] - \mathbf{E}_y[\gamma_z(\tilde{Z})]| = |\mathbf{E}_y[VA_1]|,$$

where $V \equiv \gamma_z(\tilde{Z}) - \mathbf{E}_y[\gamma_z(\tilde{Z})]$ because $\mathbf{E}_y[A_1] = 1$. The term (II) is bounded by the absolute value of

$$\int_{\bar{\lambda}_1 - 3\delta}^{\bar{\lambda}_1 + 3\delta} \mathbf{E}[VA_1|\tilde{Y} = y, \lambda_1(X) = \bar{\lambda}, \lambda_0(X) = \bar{\lambda}_0]\, dF_{\lambda_1}(\bar{\lambda}|y, \bar{\lambda}_0)$$

$$= \int_{\bar{\lambda}_1 - 3\delta}^{\bar{\lambda}_1 + 3\delta} \mathbf{E}[V|\tilde{Y} = y, \lambda_1(X) = \bar{\lambda}, \lambda_0(X) = \bar{\lambda}_0]\, dF_{\lambda_1}(\bar{\lambda}|y, \bar{\lambda}_0),$$

or by $3\bar{C}\delta$, similarly as before. This implies that (II) $\leq 3\bar{C}\delta$. $\quad\square$

**A.2. Proofs of the main results.**

PROOF OF THEOREM 1. (i) We first prove the following three claims:

C1. $\sup_{(\theta,x)\in B(\theta_0,Cn^{-1/2})\times \mathbf{R}^{d_X}} |F_{n,\theta,i}(\lambda_\theta(x)) - F_0(\lambda_0(x))| = O_P(n^{-1/2})$.

C2. $\frac{1}{\sqrt{n}}\sum_{i=1}^n \beta_u(\hat{U}_i)\{\gamma_z(\tilde{Z}_i) - \gamma_z(\tilde{Z}_i)\}\gamma_y^\perp(\tilde{Y}_i) = o_P(1)$.

C3. $\frac{1}{\sqrt{n}}\sum_{i=1}^n \beta_u(\hat{U}_i)\{\gamma_y(\tilde{Y}_i) - \gamma_y(\tilde{Y}_i)\}\{\gamma_z(\tilde{Z}_i) - \gamma_z(\tilde{Z}_i)\} = o_P(1)$.

The $o_P(1)$'s in C2 and C3 are uniform in $(u,y,z) \in [0,1]^3$.

PROOF OF C1. Let $\Lambda_n \equiv \{\lambda_\theta : \theta \in B(\theta_0, Cn^{-1/2})\}$. By Assumption 1(ii)(a),

$$(30) \quad \log N_{[\cdot]}(\varepsilon, \Lambda_n, \|\cdot\|_\infty) \leq \log N_{[\cdot]}(C\varepsilon, B(\theta_0, Cn^{-1/2}), \|\cdot\|) \leq -C\log\varepsilon,$$

because $B(\theta_0, Cn^{-1/2})$ is compact. Hence C1 follows by Lemma A3.

For the proofs for C2 and C3, we assume that $\hat{U}_i$, $\hat{Z}_i$, and $\hat{Y}_i$ are estimators using the whole sample, not leave-one-out estimators. The discrepancy due to this assumption is easily shown to be asymptotically negligible.

PROOF OF C2. Observe that

$$\Delta_{n,i}^{\mathrm{I}}(x; \hat{\theta}) \equiv \sup_{y\in\mathbf{R}} |\hat{F}_{Y|U,i}(y|F_{n,\hat{\theta},i}(\lambda_{\hat{\theta}}(x))) - F_{Y|U}(y|F_0(\lambda_0(x)))|$$

$$\leq \sup_{y\in\mathbf{R}} |\hat{F}_{Y|U,i}(y|\hat{U}_i) - F_{Y|U}(y|\hat{U}_i)|$$

$$\quad + \sup_{y\in\mathbf{R}} |F_{Y|U}(y|\hat{U}_i) - F_{Y|U}(y|U_i)|$$

$$\leq \sup_{y\in\mathbf{R}} |\hat{F}_{Y|U,i}(y|\hat{U}_i) - F_{Y|U}(y|\hat{U}_i)| + O_P(n^{-1/2})$$

by C1 and Assumption 2(ii). Take large $M > 0$ and let

$$(31) \quad v_n(u) \equiv M(h^2 1\{|u - 1| \geq h/2\} + h 1\{|u - 1| < 2h\} + n^{-1/2}h^{-1}).$$

By C1 and Lemma A4, with large probability,

$$(32) \quad \sup_{\theta \in B(\theta_0, Mn^{-1/2})} \Delta_{n,i}^{\mathrm{I}}(X_i, \theta) \leq v_n(U_i) \qquad \text{for all } i = 1, \ldots, n.$$

Take large $M$ such that $\sqrt{n}h^M = o(1)$ and expand the sum in C2 as, for example,

$$\sum_{s=1}^{M-1} B_{1n,s} + \frac{z^{M+1}}{M!\sqrt{n}} \sum_{i=1}^n \beta_u(\hat{U}_i) \exp(\bar{Z}_i z)\{\hat{Z}_i - \tilde{Z}_i\}^M \gamma_y^\perp(\tilde{Y}_i)$$
$$\equiv B_{1n} + B_{2n},$$

where $\bar{Z}_i$ lies between $\hat{Z}_i$ and $\tilde{Z}_i$ and

$$B_{1n,s} \equiv \frac{z^{s+1}}{s!\sqrt{n}} \sum_{i=1}^n \beta_u(\hat{U}_i) \exp(\tilde{Z}_i z)\{\hat{Z}_i - \tilde{Z}_i\}^s \gamma_y^\perp(\tilde{Y}_i).$$

Note that $B_{2n} = O_P(\sqrt{n}h^M) = o_P(1)$ by Lemma A4. We consider $B_{1n,s}$. Fix $c > 0$. Given $x_{(n)} \in \mathbf{R}^{nd_X}$, let $G_\theta(\cdot; x_{(n)}) \equiv G_{n,\lambda}(\lambda(\cdot); x_{(n)})$ with $\lambda = \lambda_\theta$ in (16). Denote

$$\mathcal{G}_{1n} \equiv \{G_\theta(\cdot; x_{(n)}): (x_{(n)}, \theta) \in \mathbf{R}^{nd_X} \times B(\theta_0, Mn^{-1/2})\}.$$

Then let $\mathcal{G}_n$ be the collection of $G$'s in $\mathcal{G}_{1n}$ such that $\|G - G_0\|_\infty \leq Mn^{-1/2}$, where $G_0 = F_0 \circ \lambda_0$. Define $\mathcal{B}_n \equiv \{\beta_u \circ G: (u, G) \in [0, 1] \times \mathcal{G}_n\}$. By Corollary A1 and (30),

$$(33) \quad \log N_{[\cdot]}(\varepsilon, \mathcal{B}_n, L_2(P)) \leq C\varepsilon^{-1} - C\log\varepsilon.$$

Fix a small $c > 0$ and let $\mathcal{X}_n \subset \mathbf{R}^{nd_X} \times \mathcal{G}_n$ be the collection of $(x_{(n)}, G)$'s such that $G \in \mathcal{G}_n$ and

$$(34) \quad \frac{1}{n}\sum_{j=1}^n K_h(G(x_j) - u) > c > 0.$$

Define $\mathcal{J}_{1n}^{\mathcal{G}}$ to be the set,

$$\{\phi(\cdot, G(\cdot); z_{(n)}, x_{(n)}, G): (z_{(n)}, x_{(n)}, G) \in \mathbf{R}^n \times \mathcal{X}_n\},$$

where $z_{(n)} = \{z_j\}_{j=1}^n$ and

$$\phi(z, u; z_{(n)}, x_{(n)}, G) \equiv \frac{\sum_{j=1}^n \gamma_z(z_j) K_h(G(x_j) - u)}{\sum_{j=1}^n K_h(G(x_j) - u)}.$$

Then let $\mathcal{J}_n^{\mathcal{G}} \subset \mathcal{J}_{1n}^{\mathcal{G}}$ be the class of functions $\phi(\cdot, G(\cdot))$ satisfying

$$(35) \qquad |\phi(\cdot, G(\cdot)) - F_{Z|U}(\cdot|G_0(\cdot))| \leq v_n(G_0(\cdot)).$$

Fix any arbitrarily small $w > 0$. By Assumption 3(i), Lemma A4, the fact that $|\xi_{1n}(u)| \geq 1/2$ there, and (32), from sufficiently large $n$ on,

$$P\{\hat{F}_{Z|U,i}(\cdot|F_{n,\hat{\theta}}(\lambda_{\hat{\theta}}(\cdot))) \in \mathcal{J}_n^{\mathcal{G}}\} > 1 - w$$

by increasing $M$ in the definition of $\mathcal{J}_n^{\mathcal{G}}$ sufficiently. We will now compute a bracketing entropy bound for $\mathcal{J}_n^{\mathcal{G}}$.

We apply Corollary A1 and (30) to obtain that

$$(36) \qquad \log N_{[\cdot]}(\varepsilon, \mathcal{G}_n, L_q(P)) \leq -C \log \varepsilon + C_1/\varepsilon.$$

Note that $\phi \in \mathcal{J}_n^{\mathcal{G}}$ is uniformly bounded [under the restriction (34)] and $\phi(z, u; z_{(n)}, x_{(n)}, G)$ is increasing in $z$ for all $u$ and that by (34),

$$|\phi(z, v + \eta; z_{(n)}, x_{(n)}, G) - \phi(z, v - \eta; z_{(n)}, x_{(n)}, G)| \leq C\eta/h^3.$$

Take $q > 3$ and apply Lemma A2 and (36) to obtain that

$$(37) \qquad \begin{aligned} &\log N_{[\cdot]}(\varepsilon, \mathcal{J}_n^{\mathcal{G}}, L_2(P)) \\ &\leq \log N_{[\cdot]}(C\varepsilon, \mathcal{J}_n^{\mathcal{G}}, L_q(P)) \\ &\leq -C \log \varepsilon - C \log h + C h^{-3/q} \varepsilon^{-(q+1)/q}. \end{aligned}$$

Let $\varphi_0(Z, X) \equiv F_{Z|U}(Z|G_0(X))$ and

$$f_{\beta,z,y,\varphi}(X, Z, \tilde{Y}) \equiv z^{s+1} \beta(X) e^{\varphi_0(Z,X)z} \{\varphi(Z, X) - \varphi_0(Z, X)\}^s \gamma_y^{\perp}(\tilde{Y}).$$

Let $\Pi_n \equiv \{f_{\beta,z,y,\varphi} : (\beta, z, y, \phi) \in \mathcal{B}_n \times [0,1]^2 \times \mathcal{J}_n^{\mathcal{G}}\}$. We take its envelope to be $v_n^s(G_0(\cdot))$, for which we observe that

$$(38) \qquad \mathbf{E}[v_n^{2s}(G_0(X))] = \int_0^1 v_n^{2s}(u)\, du = O(h^{2s+1}).$$

For $\Gamma \equiv \{\gamma_y : y \in [0,1]\}$, $\log N_{[\cdot]}(\varepsilon, \Gamma, L_2(du)) \leq -C \log \varepsilon$, with $du$ denoting the Lebesgue measure on $[0,1]$. Using this, (33) and (37), we conclude that

$$(39) \qquad \log N_{[\cdot]}(\varepsilon, \Pi_n, L_2(P)) \leq C h^{-3/q} \varepsilon^{-(q+1)/q} - C \log h - C \log \varepsilon.$$

Now, let us prove that $B_{1n,s} = o_P(1), s = 1, \ldots, M - 1$. With large probability,

$$B_{1n,s} \leq \sup_{f \in \Pi_n} |\sqrt{n}(\mathbf{P}_n - \mathbf{P})f| + \sup_{f \in \Pi_n} |\sqrt{n}\mathbf{P}f|,$$

where $\mathbf{P}_n$ denotes the empirical measure of $(Z_i, \tilde{Y}_i, X_i)_{i=1}^n$. Hence by the maximal inequality [Pollard (1989); see also Theorem A.2 of van der Vaart (1996)] and by (38) and (39),

$$\sup_{f \in \Pi_n} |\sqrt{n}(\mathbf{P}_n - \mathbf{P})f| \leq Ch^{-3/(2q)} \int_0^{Ch^{s+1/2}} \varepsilon^{-(q+1)/(2q)} \, d\varepsilon + o(1) = o(1).$$

We conclude that $B_{1n,s} \leq \sqrt{n} \sup_{f \in \Pi_n} |\mathbf{P}f| + o_P(1)$. By (35), the absolute value on the right-hand side is bounded by

$$(40) \qquad \mathbf{E}[|\bar{\gamma}_{z,v_n}(\tilde{Z}_i; U_i)||\mathbf{E}[\gamma_y^\perp(\tilde{Y}_i)|\tilde{Z}_i, \lambda_\theta(X_i), \lambda_0(X_i)]|] \leq \mathrm{I} + \mathrm{II},$$

where $\bar{\gamma}_{z,v_n}(\tilde{z}; u) \equiv z^{s+1} e^{\bar{z} \times z} v_n(u)^s$,

$$\mathrm{I} \equiv \mathbf{E}[|\bar{\gamma}_{z,v_n}(\tilde{Z}_i; U_i)||\mathbf{E}[\gamma_y(\tilde{Y}_i)|\tilde{Z}_i, \lambda_0(X_i)]|]$$

and

$$\mathrm{II} \equiv \mathbf{E}[|\bar{\gamma}_{z,v_n}(\tilde{Z}_i; U_i)||\mathbf{E}[\gamma_y(\tilde{Y}_i)|\tilde{Z}_i, \lambda_\theta(X_i), \lambda_0(X_i)] - \mathbf{E}[\gamma_y(\tilde{Y}_i)|\tilde{Z}_i, \lambda_0(X_i)]|].$$

Term I is bounded by

$$o(1) \times \left(\sup_{z,y,u} |\mathbf{E}[\gamma_y^\perp(\tilde{Y}_i)|\tilde{Z}_i = z, U_i = u]|\right) = o(1) \times O(n^{-1/2}) = o(n^{-1/2}),$$

under the null hypothesis or the local alternatives of the theorem. By Assumption 1(iii) with the aid of Lemma A5, term II is bounded by

$$\mathbf{E}[|\bar{\gamma}_{z,v_n}(\tilde{Z}_i; U_i)|] \times C\|\lambda_\theta - \lambda_0\|_\infty = o(n^{-1/2})$$

uniformly over $(z, y, \theta) \in [0,1]^2 \times B(\theta_0, Mn^{-1/2})$. Therefore, $B_{1n,s} = o_P(1)$, completing the proof.   $\square$

PROOF OF C3.    With large probability, the absolute value of the sum in C3 is bounded by

$$\frac{C}{\sqrt{n}} \sum_{i=1}^n \gamma_{y,v_n}(\tilde{Y}_i, U_i)\gamma_{z,v_n}(\tilde{Z}_i, U_i) + o_P(1),$$

where $\gamma_{y,v_n}(\bar{y}; u) \equiv \gamma_y(\bar{y} - v_n(u)) - \gamma_y(\bar{y} + v_n(u))$. Define $\mathcal{J}_{2n} \equiv \{\gamma_{y,v_n}(\cdot, \cdot) \times \gamma_{z,v_n}(\cdot, \cdot) : (z, y) \in [0,1]^2\}$. Note that

$$\mathbf{P}\phi^2 \leq C\mathbf{E}[\gamma_{y,v_n}(\tilde{Y}_i, U_i)\gamma_{z,v_n}(\tilde{Z}_i, U_i)] = O(h^3) \to 0 \qquad \text{for each } \phi \in \mathcal{J}_{2n}.$$

This implies that the finite-dimensional distributions of $\{\sqrt{n}(\mathbf{P}_n - \mathbf{P})\phi : \phi \in \mathcal{J}_{2n}\}$ converge in distribution to zero. By Theorem 3.3 of Ossiander (1987), $\sqrt{n}(\mathbf{P}_n - \mathbf{P})\phi$ is asymptotically tight in $l_\infty(\mathcal{J}_{2n})$ because $\phi \in \mathcal{J}_{2n}$ is uniformly bounded and

$$(41) \qquad \log N_{[\cdot]}(\varepsilon, \mathcal{J}_{2n}, L_2(P)) \leq 2\log N_{[\cdot]}(C\varepsilon, \Gamma, L_2(P)) \leq -C\log\varepsilon.$$

Hence $\sup_{\phi \in \mathcal{J}_{2n}} |\sqrt{n}(\mathbf{P}_n - \mathbf{P})\phi| = o_P(1)$. We are left with $\sqrt{n}\mathbf{P}\phi$ to analyze. Let $1_n = 1\{|U_i - 1| \geq 2h\}$. Under the null hypothesis or the local alternatives in the theorem,

$$
\begin{aligned}
\sqrt{n}\,\mathbf{P}\phi &= \sqrt{n}\,\mathbf{E}[\gamma_{y,v_n}(\tilde{Y}_i, U_i)\gamma_{z,v_n}(\tilde{Z}_i, U_i)] \\
&= \sqrt{n}\,\mathbf{E}[\gamma_{y,v_n}(\tilde{Y}_i, U_i)\mathbf{E}[\gamma_{z,v_n}(\tilde{Z}_i, U_i)|\tilde{Y}_i, U_i]1_n] \\
&\quad + \sqrt{n}\,\mathbf{E}[\gamma_{y,v_n}(\tilde{Y}_i, U_i)\mathbf{E}[\gamma_{z,v_n}(\tilde{Z}_i, U_i)|\tilde{Y}_i, U_i](1 - 1_n)] \\
&= O(\sqrt{n}w_n^2) + O(\sqrt{n}h^3) + O(h) = o(1), \qquad w_n \equiv n^{-1/2}h^{-1},
\end{aligned}
$$

because the expectation of $1\{|U_i - 1| < 2h\}$ is $O(h)$ and $n^{-1/2}h^{-2} + n^{1/2}h^3 \to 0$ as $n \to \infty$ by Assumption 3(ii)(b). Hence $\sup_{\phi \in \mathcal{J}_{2n}} \sqrt{n}\,\mathbf{P}\phi = o(1)$, establishing C3.

Now we turn to the proof of Theorem 1. Without loss of generality, let $\beta_u$ be monotone decreasing. First we write $\hat{\nu}_n(r) - \nu_n(r)$ as

$$
\begin{aligned}
&\frac{1}{\sqrt{n}}\sum_{i=1}^{n}\beta_u(\hat{U}_i)\{\gamma_z^{\perp}(\hat{Z}_i)\gamma_y^{\perp}(\hat{Y}_i) - \gamma_z^{\perp}(\tilde{Z}_i)\gamma_y^{\perp}(\tilde{Y}_i)\} \\
&\quad + \frac{1}{\sqrt{n}}\sum_{i=1}^{n}\{\beta_u(\hat{U}_i) - \beta_u(U_i)\}\gamma_z^{\perp}(\tilde{Z}_i)\gamma_y^{\perp}(\tilde{Y}_i) \\
&\equiv A_{1n} + A_{2n}.
\end{aligned}
\tag{42}
$$

We show that $A_{jn} = o_P(1)$, $j = 1, 2$. Write $\gamma_z^{\perp}(\hat{Z}_i)\gamma_y^{\perp}(\hat{Y}_i) - \gamma_z^{\perp}(\tilde{Z}_i)\gamma_y^{\perp}(\tilde{Y}_i)$ as

$$
\begin{aligned}
&\{\gamma_z^{\perp}(\hat{Z}_i) - \gamma_z^{\perp}(\tilde{Z}_i)\}\gamma_y^{\perp}(\tilde{Y}_i) \\
&\quad + \{\gamma_z^{\perp}(\hat{Z}_i) - \gamma_z^{\perp}(\tilde{Z}_i)\}\{\gamma_y^{\perp}(\hat{Y}_i) - \gamma_y^{\perp}(\tilde{Y}_i)\} \\
&\quad + \{\gamma_y^{\perp}(\hat{Y}_i) - \gamma_y^{\perp}(\tilde{Y}_i)\}\gamma_z^{\perp}(\tilde{Z}_i).
\end{aligned}
\tag{43}
$$

By decomposing $A_{1n}$ into three terms according to the decomposition in (43) and applying C2 and C3, we obtain that $A_{1n} = o_P(1)$.

As for $A_{2n}$, define

$$
\mathcal{J}_{1n} \equiv \{((\beta_u \circ G) - \beta_u)\gamma_y^{\perp}\gamma_z^{\perp} : (G, u, y, z) \in \mathcal{G}_n \times [0, 1]^3\}.
$$

Following the arguments used to show $\sup_{f \in \Pi_n} |\sqrt{n}(\mathbf{P}_n - \mathbf{P})f| = o_P(1)$ in the proof of C2, we can write

$$
|A_{2n}| \leq \sup_{f \in \mathcal{J}_{1n}} \sqrt{n}|\mathbf{P}f| + o_P(1).
$$

Similarly, as in the steps in and below (40), the leading term above is $o_P(1)$ due to the fact that

$$
\sup_{(z,y,u) \in [0,1]^3} |\mathbf{E}[\gamma_y^{\perp}(\tilde{Y}_i)\gamma_z^{\perp}(\tilde{Z}_i)|U_i = u]| = O(n^{-1/2})
$$

under the null hypothesis and the local alternatives and the fact that

$$\sup_{(G,u)\in\mathcal{G}_n\times[0,1]}\mathbf{E}[|\beta_u(G(X_i))-\beta_u(U_i)|]$$

$$\leq\sup_{u\in[0,1]}\int_0^1\{\beta_u(\bar{u}+Mn^{-1/2})-\beta_u(\bar{u}-Mn^{-1/2})\}\,d\bar{u}=O(n^{-1/2}),$$

which follows because $\beta_u$ is increasing and bounded in $[0,1]$.

(ii) Clearly the class $\mathcal{J}\equiv\{\beta_u\gamma_y^\perp\gamma_z^\perp:(u,y,z)\in[0,1]^3\}$ is $P$-Donsker because its element is a product of uniformly bounded monotone functions. The weak convergence of $\nu_n(r)$ immediately follows, and the weak convergence of $\hat{\nu}_n(r)$ follows from (i). $\square$

PROOF OF THEOREM 2. Let $\nu_{n,b}^{0*}(r)\equiv\frac{1}{\sqrt{n}}\sum_{i=1}^n\omega_{i,b}\beta_u(U_i)\gamma_z^\perp(\tilde{Z}_i)\gamma_y^\perp(\tilde{Y}_i)$. Also let $P_\omega$ be the distribution of $(\omega_{i,b})_{i=1}^n$ and $\mathbf{E}_\omega$ be the associated expectation. It suffices to show that for each $\varepsilon>0$,

$$(44)\qquad P_\omega\Big\{\sup_{r\in[0,1]^3}|\nu_{n,b}^*(r)-\nu_{n,b}^{0*}(r)|>\varepsilon\Big\}=o_P(1).$$

This is because the class $\mathcal{J}$ is $P$-Donsker as we saw in the proof of Theorem 1(ii), and by the conditional multiplier uniform CLT of Ledoux and Talagrand (1988) [e.g., Theorem 2.9.7 in van der Vaart and Wellner (1996)],

$$d(F_{\Gamma\nu_n^{0*}}^*,F_{\Gamma\nu})=o(1)\qquad\text{a.s.}$$

We turn to (44). Let $\hat{S}_i=(U_i,Z_i,Y_i,\hat{U}_i,\hat{Z}_i,\hat{Y}_i)$ and

$$\xi_i(r)\equiv\xi(\hat{S}_i;r)\equiv\beta_u(U_i)\gamma_z^\perp(\tilde{Z}_i)\gamma_y^\perp(\tilde{Y}_i)-\beta_u(\hat{U}_i)\gamma_z^\perp(\hat{Z}_i)\gamma_y^\perp(\hat{Y}_i).$$

Then note that for all $r\in[0,1]^3$,

$$(45)\qquad\mathbf{E}_\omega\bigg[\bigg|\frac{1}{\sqrt{n}}\sum_{i=1}^n\omega_{i,b}\xi(\hat{S}_i;r)\bigg|^2\bigg]=\frac{1}{n}\sum_{i=1}^n\xi(\hat{S}_i;r)^2\to_P 0,$$

using the proof of Theorem 1. Let $\rho_n(r_1,r_2)=\sqrt{\frac{1}{n}\sum_{i=1}^n(\xi_i(r_1)-\xi_i(r_2))^2}$. Observe that

$$-\frac{\sqrt{5}-1}{2}|\xi_i(r_1)-\xi_i(r_2)|\leq\omega_{i,b}(\xi_i(r_1)-\xi_i(r_2))\leq\frac{\sqrt{5}+1}{2}|\xi_i(r_1)-\xi_i(r_2)|.$$

By Hoeffding's inequality [e.g., Lemma 3.5 of van de Geer (2000), page 33],

$$P_\omega\bigg\{\bigg|\frac{1}{\sqrt{n}}\sum_{i=1}^n\omega_{i,b}(\xi_i(r_1)-\xi_i(r_2))\bigg|>\varepsilon\bigg\}\leq\exp\bigg(-\frac{4\varepsilon^2}{(\sqrt{5}+1)^2\rho_n(r_1,r_2)^2}\bigg).$$

Therefore, the process $\nu_{n,b}^*(r)$ is sub-Gaussian with respect to $\rho_n$.

Let us compute the covering number of $[0,1]^3$ with respect to $\rho_n$. Let $P_n$ be the empirical measure of $(U_i, Z_i, Y_i, \hat{U}_i, \hat{Z}_i, \hat{Y}_i)_{i=1}^n$. Observe that $\xi(\cdot; r)$ is bounded variation, so that for $\Xi = \{\xi(\cdot; r) : r \in [0,1]^3\}$,

$$\log N_{[\cdot]}(\varepsilon, \Xi, L_2(P_n)) \leq \frac{C}{\varepsilon}.$$

Since $\rho_n^2(r, r_j) = \|\xi(\cdot; r) - \xi(\cdot; r_j)\|_{P_n, 2}^2$, this implies that $\log N(\varepsilon, [0,1]^3, \rho_n) \leq C/\varepsilon$.

Now, using Corollary 2.2.8 of van der Vaart and Wellner [(1996), page 101], for any $r_0 \in [0,1]^3$,

$$
\begin{aligned}
(46) \quad & \mathbf{E}_\omega \left[ \sup_{r \in [0,1]^3} \left| \frac{1}{\sqrt{n}} \sum_{i=1}^n \omega_{i,b} \xi(\hat{S}_i; r) \right| \right] \\
& \leq \mathbf{E}_\omega \left[ \left| \frac{1}{\sqrt{n}} \sum_{i=1}^n \omega_{i,b} \xi(\hat{S}_i; r_0) \right| \right] + C \int_0^\infty \sqrt{\log D(\varepsilon, \rho_n)} \, d\varepsilon,
\end{aligned}
$$

where $C$ is an absolute constant and $D(\varepsilon, \rho_n)$ is an $\varepsilon$-packing number of $[0,1]^3$ with respect to $\rho_n$. The leading term on the right-hand side vanishes in probability by (45). As for the second term, note that

$$\sup_{r_1, r_2 \in [0,1]^3} \rho_n(r_1, r_2) \to_P 0 \qquad \text{as } n \to \infty$$

from the proof of Theorem 1. Therefore, we can take $\delta_n \to 0$ such that

$$P\left\{ \sup_{r_1, r_2 \in [0,1]^3} \rho_n(r_1, r_2) < \delta_n \right\} \to 1.$$

With probability approaching one,

$$\int_0^\infty \sqrt{\log D(\varepsilon, \rho_n)} \, d\varepsilon = \int_0^{C\delta_n} \sqrt{\log D(\varepsilon, \rho_n)} \, d\varepsilon \leq C \int_0^{C\delta_n} \varepsilon^{-1/2} \, d\varepsilon \to 0.$$

We conclude that the expectation on the left-hand side of (46) vanishes in probability. We obtain (44). $\square$

PROOF OF THEOREM 3. We first show that $\sup_{r \in [0,1]^2 \times \mathcal{Z}} |\bar{\nu}_n(r) - \bar{\nu}_n^0(r)| = o_P(1)$ both under $H_0$ and the local alternatives in the theorem, where

$$\bar{\nu}_n^0(r) = \frac{1}{\sqrt{n}} \sum_{i=1}^n \frac{\beta_u(U_i)\{1\{Z_i = z\} - p_z(U_i)\}\gamma_y^\perp(\tilde{Y}_i)}{\sqrt{p_z(U_i) - p_z(U_i)^2}}.$$

By Lemma A4, $\sup_{u \in [0,1]} \sup_{z \in \mathcal{Z}} |\hat{p}_{z,i}(u) - p_z(u)| \leq v_n(u)$ where $v_n(u)$ is as defined in (31). By Assumption 1D(iv), $\hat{p}_{z,i}(u) \in (\varepsilon/2, 1 - \varepsilon/2), \varepsilon > 0$, with probability approaching one. We introduce $G_\theta(\cdot)$ and $\mathcal{G}_n$ as in the proof of Theorem 1. Define

$$\phi_z(u) \equiv \frac{\sum_{j=1}^n 1\{z_j = z\} K_h(G_\theta(x_j) - u)}{\sum_{j=1}^n K_h(G_\theta(x_j) - u)}.$$

Let $\mathcal{J}_{n,2}^{\mathcal{G}}$ be the class of functions $\phi_z(G_\theta(\cdot))$ with $(\{(x_j^\top, z_j)^\top\}_{j=1}^n, \theta)$ running in $(\mathbf{R}^{d_X} \times \mathcal{Z})^n \times B(\theta_0, Mn^{-1/2})$ and such that for some $c > 0$,

$$\frac{1}{n}\sum_{j=1}^n K_h(G_\theta(x_j) - u) > c > 0$$

and $|\phi_z(G_\theta(\cdot)) - p_z(\cdot)| \le v_n(G_0(\cdot))$. By Lemma A2, for any $w > 0$, we can find large $M$ such that

$$P\{\hat{p}_{z,i}(F_{n,\hat{\theta},i}(\lambda_{\hat\theta}(\cdot))) \in \mathcal{J}_{n,2}^{\mathcal{G}}\} > 1 - w,$$

$\log N_{[\,]}(\varepsilon, \mathcal{J}_{n,2}^{\mathcal{G}}, L_2(P)) \le Ch^{-3/q}\varepsilon^{-(q+1)/q}$ for some $q > 3$ and for all $\pi(\cdot) \in \mathcal{J}_{n,2}^{\mathcal{G}}$,

$$0 < \varepsilon/2 < \inf_{x \in \mathbf{R}^{d_X}} |\pi(x)| \le \sup_{x \in \mathbf{R}^{d_X}} |\pi(x)| < 1 - \varepsilon/2 < 1$$

from some sufficiently large $n$ on. Define

$$\acute{\mathcal{B}_n} \equiv \{\bar\beta_{\bar u, \bar z}(G(\cdot), \cdot, \cdot; \pi) : (\bar u, \bar z, \pi, G) \in [0,1] \times \mathcal{Z} \times \mathcal{J}_{n,2}^{\mathcal{G}} \times \mathcal{G}_n\},$$

where $\bar\beta_{\bar u, \bar z}(u, z, x; \pi) \equiv \beta_{\bar u}(u)\{1\{z = \bar z\} - \pi(x)\}/\sqrt{\pi(x) - \pi^2(x)}$. Then we can write

$$\bar\nu_n(r) = \frac{1}{\sqrt{n}}\sum_{i=1}^n \bar\beta_{u,z}(\hat{U}_i, Z_i, X_i; \hat{p}_{z,i} \circ (F_{n,\hat\theta,i} \circ \lambda_{\hat\theta}))\gamma_y^\perp(\hat{Y}_i).$$

To show that it is asymptotically equivalent to $\bar\nu_n^0(r)$ uniformly over $r \in [0,1]^2 \times \mathcal{Z}$, we consider the following process:

$$(47) \qquad \frac{1}{\sqrt{n}}\sum_{i=1}^n \beta(Z_i, X_i)\gamma_y^\perp(\hat{Y}_i), \qquad (y, \beta) \in [0,1] \times \acute{\mathcal{B}_n}.$$

Note that $\bar\beta_{\bar u, \bar z}(u, z, x; \pi)$ is Lipschitz in $\pi \in (\varepsilon/2, 1 - \varepsilon/2)$ with a uniformly bounded coefficient. Therefore, $\log N_{[\,]}(\varepsilon, \acute{\mathcal{B}_n}, L_2(P))$ is bounded by

$$(48) \qquad \log N_{[\,]}(C\varepsilon, \mathcal{J}_{n,2}^{\mathcal{G}}, L_2(P)) + \log N_{[\,]}(C\varepsilon, \mathcal{B}_n, L_2(P))$$
$$\le Ch^{-3/q}\varepsilon^{-(q+1)/q} + C\varepsilon^{-1} - C\log\varepsilon,$$

where $\mathcal{B}_n$ was defined prior to (33). Using this bound and proceeding with the proof of Theorem 1(i) by replacing $\gamma_{\bar z}^\perp$ by 1 and $\mathcal{B}_n$ by $\acute{\mathcal{B}_n}$, we conclude that $\bar\nu_n$ is asymptotically equivalent to $\bar\nu_n^0$.

Let $\mathcal{J}_z \equiv \{\beta_u(\cdot)\gamma_{\bar{y}}^{\perp}(\cdot)\{1\{\cdot = z\} - p_z(\cdot)\}/\sqrt{p_z(\cdot) - p_z(\cdot)^2} : (y, u) \in [0, 1]^2\}$. Then,

$$\log N_{[\cdot]}(\varepsilon, \mathcal{J}_z, L_2(P)) \le C\varepsilon^{-1}.$$

Therefore, $\bigcup_{z \in \mathcal{Z}} \mathcal{J}_z$ is $P$-Donsker because $\mathcal{Z}$ is a finite set. The weak convergence of $\bar{\nu}_n^0$ to $\bar{\nu}$ follows. The weak convergence of $\bar{\nu}_n$ to $\bar{\nu}$ follows from its asymptotic equivalence with $\bar{\nu}_n^0$. $\square$

PROOF OF THEOREM 4. The proof is almost the same as the proof of Theorem 2. We omit the details. $\square$

**Acknowledgments.** I thank an Associated Editor and a referee for their comments. I am grateful to Yoon-Jae Whang and Guido Kuersteiner who gave me valuable advice at the initial stage of this research. All errors are mine.

DEPARTMENT OF ECONOMICS
UNIVERSITY OF PENNSYLVANIA
528 MCNEIL BUILING
3718 LOCUST WALK
PHILADELPHIA, PENNSYLVANIA 19104
USA
E-MAIL: kysong@sas.upenn.edu